\documentclass[12pt]{article}
\usepackage{rotating}
\usepackage{color}
\usepackage{amssymb}
\usepackage{amsmath,amsthm}
\usepackage{amsfonts}
\usepackage{graphicx}
\usepackage{graphics}
\usepackage{setspace}
\numberwithin{equation}{section}

\begin{document}

\date{\small\textsl{\today}}
\title{The BEM and DRBEM schemes for the numerical solution of the two-dimensional time-fractional diffusion-wave equations}
\author{ \large Peyman Alipour \footnote{Corresponding author.\newline {\em  E-mail
addresses:} palipour@stevens.edu. \newline} $\vspace{.2cm} $ \large \\
\small{\em School of Business, Stevens Institute of Technology, Hoboken, NJ, USA}\vspace{-1mm}} \maketitle
\vspace{.9cm}
\begin{abstract}
In this paper we apply the boundary elements method (BEM) and the dual reciprocity  boundary elements method (DRBEM) for the numerical
solution of two--dimensional time--fractional partial differential equations (TFPDEs). The fractional derivative of problem is described in the Caputo sense.
In BEM, the main equation deduces to solving the Helmholtz equation in each time step. Therefore, we should compute the domain integral in each time step. So, we presented an approach to compute the domain integral with no singularity.
On the other hand the DRBEM has the flexibility of discretizing only the boundary of the computational domain and evaluates the solution at any required interior point.
We employ the radial basis functions (RBFs) for interpolation of the inhomogeneous and time derivative terms.
The proposed method is employed for solving some problems in two--dimensions on unit square and some other complex regions to demonstrate the efficiency of the proposed method. \\
\vspace{.5cm}

\end{abstract}

\section{Introduction}
All natural phenomenons can not be modeled by partial differential equations of integer-order, therefore some of them modeled by fractional partial differential equations (FPDEs) \cite{A. V,R. Metzler,I. Padlubny}. A history of the development of fractional differential operators can be found in \cite{Miller,Oldham}. Oldham \cite{Oldham} and Podlubny \cite{I. Padlubny} have pointed to some applications of fractional partial differential equations. The analytical results and detailed physical description for the fractional diffusion-wave problem can be found in \cite{R. Gorenflo,F. Liu1, F. Mainardi3, W. R. Schneider, W. Wyss}.
Analytical solutions for most fractional differential equations are not usually expressed explicitly, therefore many researchers have utilized numerical solution strategies based on convergence and stability analysis \cite{K. Diethelm,Liu,F. Liu,Tadjeran,S.B. Yuste,Zhuang}.
Zhang \cite{Yang Zhang} used an implicit unconditional stable difference scheme for a kind of linear space-time fractional convection-diffusion equation.
Meerschaert and Tadjeran \cite{Mark M. Meerschaert} examined some practical numerical methods to solve a class of initial-boundary value FPDEs with variable coefficients on a finite domain.
Jiang and Ma \cite{Y. Jiang} developed high-order methods for solving time-fractional partial differential equations (TFPDEs). The proposed high-order method is based on high-order finite element method for space direction and finite difference method for time variable.
Cui \cite{Cui1} solved the one-dimensional fractional diffusion equation via high-order compact finite difference scheme and obtained a fully discrete implicit scheme by Grunwald-Letnikov discretization of the Riemann-Liouville derivative. Author of \cite{Cui2} also used a high-order compact exponential scheme for solving the fractional convection-diffusion equation.
Mohebi et al. \cite{A. Mohebbi} utilizied a meshless technique based on collocation and radial basis functions (RBFs) for solving the fractional nonlinear Schrodinger equation arising in quantum mechanics. They also used a high-order compact finite difference method for solving TFPDEs \cite{A. Mohebbi1, A. Mohebbi2}. \\
The boundary elements method (BEM) \cite{Ang1,C. A. Brebbia,J. Katsikadelis} is a suitable numerical method for solving boundary value problems.
The basic idea of this method is the transformation of the original differential equation describing the behaviour of the unknown both inside and on the boundary of the domain into an equivalent integral equation only on the boundary. BEM requires only the boundary discretization which reduces the dimension of the problem under consideration by one. Thus, a smaller system of equations is obtained in comparison with the numerical methods requiring domain discretization.
The method has become a widely used technique with the applications around in many areas of engineering such as fluid mechanics, magnetohydrodynamic, electrodynamics and so on, for instance see \cite{N. Alsoy-Akun,Ang,C. Bozkaya,C. Bozkaya1,H. Hosseinzadeh,D. Mirzaei,M. Tezer--Sezgin,M. Tezer--Sezgin1,Z. Sedaghatjoo1}.
The BEM always requires a fundamental solution to the original differential equation in order to avoid domain integrals in the formulation
of the boundary integral equation also the nonhomogeneous and nonlinear terms are incorporated in the formulation by means of domain integrals.
The dual reciprocity boundary elements method (DRBEM) has the flexibility of discretizing only the boundary of the computational domain and evaluates the solution at any required interior point.
DRBEM does not need the mesh generation of the computational domain and does not require the fundamental solution of the main equation (only the fundamental solution of Laplace equation is needed). Some engineering problems are solved by the DRBEM, for instance, see \cite{Ang,Gumgum}.
The DRBEM has been developed for the numerical solution of TFPDEs by Katsikadelis.
He solved the general linear fractional diffusion--wave equation in bounded inhomogeneous anisotropic bodies \cite{Katsikadelis2}.
He also in \cite{Katsikadelis1} presented a numerical method for the solution of FPDEs. The solution procedure applies to both linear and nonlinear problems described by evolution type equations involving time--fractional derivatives in bounded domains of arbitrary shape.

Fractional diffusion-wave equations are a class of important TFPDEs which has been widely applied in the modeling of an anomalous diffusion, sub-diffusion systems, fractals and percolation clusters, polymers, random and disordered media, geophysical and geological processes \cite{ V. V. Gafiychuk, R. Metzler}.
As mentioned in \cite{H.Jiang}, the time-fractional diffusion-wave and diffusion equations can be written in the following form:
\begin{equation}\label{formula 1}
\begin{array}{l}
D_{t}^{\alpha}u(X,t)=\kappa \nabla ^{2} u(X,t) +g(X,t), \ \ \ \ X \in \Omega, \ \ t > 0,
\end{array}
\end{equation}
along with initial conditions
\begin{equation}\label{formula 2}
\begin{array}{l}
u(X,0)=\phi(X), \ \ \ \ \ \displaystyle \frac{\partial u(X,0)}{\partial t}=\psi(X), \ \ \ \ X \in \Omega,
\end{array}
\end{equation}
and the boundary conditions
\begin{equation}\label{formula 3}
\begin{array}{l}
u(X,t)=\varphi(X,t), \ \ \ \ \ X \in \Gamma,
\end{array}
\end{equation}
where $X$ and $t$ are the space and time variables, respectively, $\Gamma$ is the boundary enclosing $\Omega$, $\kappa$ is an arbitrary positive constant, $g(X,t)$ is a sufficiently smooth function, $0 < \alpha < 2$ and $D_{t}^{\alpha}$ is a Caputo fractional derivative of order $\alpha$ defined as \cite{I. Padlubny}
\begin{equation}\label{formula 4}
D_{t}^{\alpha}u(X,t)=\displaystyle \frac{\partial^{\alpha} u(X,t)}{\partial t^{\alpha}} = \left\{ {\begin{array}{*{20}{c}}
{{\displaystyle \frac{1}{\Gamma(m-\alpha)}  \displaystyle \int_{0}^{t} \frac{u^{(m)}(X,\tau)}{(t-\tau)^{1+\alpha-m}} dt,}}&{m-1 < \alpha < m},\\
\\
{{\displaystyle \frac{d^{m}}{dt^{m}} u(X,t), \ \ \ \ \ \ \ \ \ \ \ \ \ \ \ \ \ \ \ \ \ \ \ \ \ }}&{\alpha=m \in \mathbb{N}.}
\end{array}} \right.
\end{equation}
There are some numerical methods which have been applied for solving the fractional diffusion-wave equations, for instance, finite difference \cite{Ya-Nan Zhang, Zhi-zhong Sun} and finite element \cite{Limei Li} methods.
In addition, the authors in \cite{F. Mainardi4} by using Laplace transform solved fractional diffusion-wave equation by converting the main equation into a system of ordinary differential equations.

The outline of this paper is as follows. In Section 2, the time discretization scheme has been proposed. In this section, we also show that the discretization deduces the FPDEs to an inhomogeneous Helmholtz equation. In Section 3, BEM has been applied for the fractional diffusion--wave equation.
Computation of domain integral is presented in Subsection 3.1.
In Section 4, by using the DRBEM, the discretized version of the fractional diffusion--wave equation, is obtained and an iterative scheme is described for the time derivatives. The numerical results are presented in Section 5. Section 6 completes this paper with a brief conclusion.

\section{Time discrete scheme}
At first, we define the functional space endowed with standard norms and inner products.
Let $\Omega$ denote a bounded and open domain in $\mathbb{R^{d}}$, for $d=2$ or $3$. For $p < + \infty$, we denote by $L^{p}(\Omega)$ the space of measurable functions $u : \Omega \rightarrow \mathbb{R}$ such that $\int_{\Omega} | u(x) |^{p} dx \leq \infty$. It is a Banach space for the norm \cite{H. Brezis}

\begin{equation*}
\begin{array}{l}
\parallel u \parallel _{L^{p}(\Omega)}= \left ( \displaystyle  \int_{\Omega} |u(x)|^{p} dx  \right )^{\frac{1}{p}}.
\end{array}
\end{equation*}
The space $L^{2}(\Omega)$ is a Hilbert space with the inner product
\begin{equation*}
\begin{array}{l}
(u,v)=\displaystyle \int_{\Omega} u(x) v(x) dx,
\end{array}
\end{equation*}
which induces the norm
\begin{equation*}
\begin{array}{l}
\parallel u \parallel _{L^{2}(\Omega)}= \left ( \displaystyle  \int_{\Omega} |u(x)|^{2} dx  \right )^{\frac{1}{2}}.
\end{array}
\end{equation*}
And
\begin{equation*}
\begin{array}{l}
H^{1}(\Omega)= \displaystyle \{ \upsilon \in L^{2}(\Omega),  \displaystyle \frac{d\upsilon}{dx} \in L^{2}(\Omega) \displaystyle \}, \\
\\
H_{0}^{1}(\Omega)= \displaystyle \{ \upsilon \in H^{1}(\Omega), \upsilon |_{\Gamma}=0  \displaystyle \}.
\end{array}
\end{equation*}
For discretization of the time variable, we need some preliminary. We define
\begin{equation*}
\begin{array}{l}
t_{k}:=k \tau, \ \ \ \ k=0,1,2,...,N_{t},
\end{array}
\end{equation*}
where $\tau=\displaystyle \frac{T}{N_{t}}$ and
\begin{equation}\label{formula 5}
\begin{array}{l}
\xi_t u^{n-\frac{1}{2}}=\displaystyle \frac{1}{\tau}(u^n - u^{n-1}), \ \ \ \ v(X,t) = \displaystyle  \frac{\partial u(X,t)}{\partial t}.
\end{array}
\end{equation}
The difference scheme that we will consider for solving problem (\ref{formula 1})-(\ref{formula 3}) is as follows \cite{Zhi-zhong Sun}:
\begin{equation}\label{formula 6}
\begin{array}{l}
\displaystyle \frac {1}{\tau \Gamma(2-\alpha)} \left[ a_{0} \xi_{t} u^{n-\frac{1}{2}} - \displaystyle \sum_{k=1}^{n-1} (a_{n-k-1}-a_{n-k}) \xi_{t} u^{k-\frac{1}{2}} - a_{n-1}\psi \right] = \kappa \nabla^{2} u^{n} + g^{n-\frac{1}{2}}, \ \ n \geq 1,
\end{array}
\end{equation}
where
\begin{equation}\label{formula 7}
\begin{array}{l}
a_k := \displaystyle \int_{t_{k}}^{t_{k+1}} \frac{dt}{t^{\alpha-1}} = \frac{\tau^{2-\alpha}}{2-\alpha}[(k+1)^{2-\alpha}-k^{2-\alpha}], \ \  k \geq 0,
\end{array}
\end{equation}
and
\begin{equation*}
\begin{array}{l}
g^{n-\frac{1}{2}}=g(X,\frac{t_{n}+t_{n-1}}{2}), \ \ \ \ n \geq 1.
\end{array}
\end{equation*}
We discrete the time variable using the following lemmas.\\
\textbf{Lemma 1}. \cite{Zhi-zhong Sun} Suppose $g(t) \in C^{2}[0,t_{n}]$ and $1 < \alpha < 2$ then we have
\begin{equation*}
\begin{array}{l}
\left| \displaystyle \int_{\Omega} \displaystyle \frac{g^{'}(t)}{(t_{n}-t)^{\alpha-1}} dt - \displaystyle \frac{1}{\tau} \left[ a_{0} g(t_{n}) - \displaystyle \sum_{k=1}^{n-1} (a_{n-k-1}-a_{n-k})g(t_{k}) - a_{n-1} g(t_{0}) \right] \right| \\
\leq \displaystyle \frac{1}{2-\alpha} \left[ \displaystyle \frac{2-\alpha}{12} + \displaystyle \frac{2^{3-\alpha}}{3-\alpha}
- (1+2^{1-\alpha})\right] \mathop {\max  | {g''(t)}|{\tau ^{3 - \alpha
}}}\limits_{0 \le t \le {t_n}} .
\end{array}
\end{equation*}
We introduce the following notation:
\begin{equation}\label{formula 8}
w(X,t) := \frac{1}{{\Gamma (2 - \alpha )}}\int_0^{t}
\frac{\partial v(X,s)}{\partial s} \frac{ds}{(t -s)^{\alpha  - 1}},
\end{equation}
and
\begin{equation*}
\mathcal{P}  \left( u^{n-\frac{1}{2}},q \right) := a_{0} u^{n-\frac{1}{2}}  - \displaystyle \sum_{k=1}^{n-1} (a_{n-k-1}-a_{n-k}) u^{k-\frac{1}{2}} - a_{n-1} q.
\end{equation*}
So from Lemma 1 we have
\begin{equation*}
w^{n}= \displaystyle \frac{1}{\Gamma(2-\alpha)} \displaystyle \int _{0}^{t_{n}} \displaystyle \frac{\partial v(X,s)}{\partial s} \displaystyle \frac{dt}{(t_{n}-s)^{\alpha-1}} = \displaystyle \frac{1}{\tau \Gamma(2-\alpha)} \mathcal{P}(v^{n},\psi) + O(\tau^{3-\alpha}),
\end{equation*}
such that $w^{n}=w(X,t_{n})$.
Using Taylor expansion, it follows from (\ref{formula 5}) and (\ref{formula 8}) that
\begin{equation}\label{formula 9}
v^{n-\frac{1}{2}} = \xi_{t} u^{n-\frac{1}{2}} + (r_{1})^{n-\frac{1}{2}},
\end{equation}
\begin{equation}\label{formula 10}
w^{n-\frac{1}{2}} = \displaystyle \frac{1}{2} \left( w^{n} + w^{n-1}  \right) = \kappa \nabla^{2} u^{n} + g^{n-\frac{1}{2}}+(r_{2})^{n-\frac{1}{2}}, \ \ \ \ \ \ X \in \Omega, \ \ \ \ n \geq 1,
\end{equation}
where $u^{n}=u(X,t_{n})$ and there exists a constant $c_{1}$ such that
\begin{equation}\label{formula 11}
\left| (r_{1})^{n-\frac{1}{2}}  \right| \leq c_{1} \tau^{2}, \ \ \ \
\left| (r_{2})^{n-\frac{1}{2}}  \right| \leq c_{1} (\tau).
\end{equation}
Using Lemma 1 and $v^{0}=v(X,0)=\psi(X)=\psi$, the above relation can be written as
\begin{equation*}
w^{n-\frac{1}{2}} = \displaystyle \frac{1}{\tau \Gamma(2-\alpha)}  \mathcal{P}(v^{n-\frac{1}{2}},\psi) + O(\tau^{3-\alpha}).
\end{equation*}
Consequently
\begin{equation}\label{formula 12}
w^{n-\frac{1}{2}}  = \displaystyle \frac{1}{\tau \Gamma(2-\alpha)}  \mathcal{P}(v^{n-\frac{1}{2}},\psi) + (r_{3})^{n-\frac{1}{2}},
\end{equation}
where there exists a constant $c_{2}$ such that \cite{Zhi-zhong Sun}
\begin{equation}\label{formula 13}
\left| (r_{3})^{n-\frac{1}{2}} \right| \leq c_{2} \tau ^{3-\alpha}.
\end{equation}
If we substitute (\ref{formula 12}) into (\ref{formula 9}), we obtain
\begin{equation*}
w^{n-\frac{1}{2}}  = \displaystyle \frac{1}{\tau \Gamma(2-\alpha)}  \mathcal{P}(\xi_{t} u^{n-\frac{1}{2}},\psi) + \displaystyle \frac{1}{\tau \Gamma(2-\alpha)}  \mathcal{P}\left((r_{1})^{n-\frac{1}{2}},0\right) + (r_{3})^{n-\frac{1}{2}}.
\end{equation*}
Then substituting the above result into (\ref{formula 1}), yields
\begin{equation}\label{formula 14}
\begin{array}{l}
\displaystyle \frac{1}{\tau \Gamma(2-\alpha)}  \mathcal{P}(\xi_{t} u^{n-\frac{1}{2}},\psi) =\kappa \nabla^{2} u^{n} + g^{n-\frac{1}{2}} + R^{n-\frac{1}{2}},\ \ \ \  n \geq 1,
\end{array}
\end{equation}
where
\begin{equation*}
R^{n-\frac{1}{2}} =  \left\{ \displaystyle \frac{1}{\tau \Gamma(2-\alpha)}  \mathcal{P}\left((r_{1})^{n-\frac{1}{2}},0\right) + (r_{3})^{n-\frac{1}{2}} \right\} + (r_{2})^{n-\frac{1}{2}}.
\end{equation*}
Now by omitting the small term $R^{n-\frac{1}{2}}$, we can construct the following difference scheme for problem (\ref{formula 14}),
\begin{equation}\label{formula 16}
\left\{ \begin{array}{l}
\displaystyle \frac{1}{\tau \Gamma(2-\alpha)}  \mathcal{P}(\xi_{t} u^{n-\frac{1}{2}},\psi)  = \kappa \nabla^{2} u^{n} + g^{n-\frac{1}{2}} , \ \ \ \ \ X \in \Omega, \ \ \ \ \ \ n \geq 1,\\
u^{n}= \varphi (X,t_{n}),\ \ \ \ X \in \Gamma,
\end{array} \right.
\end{equation}
where $u^{n}$ is the approximation solution of (\ref{formula 14}).
So, by time discretization, solving Eq. (\ref{formula 1}) leads to solve Eq. (\ref{formula 16}). In the other hand by using the definition of  $ \xi_t u^{n-\frac{1}{2}} $, as presented in Eqs. (\ref{formula 5}), the following expansion will be obtained

\begin{equation}\label{formula 17}
\begin{array}{l}
\displaystyle \frac{1}{\tau \Gamma(2-\alpha)}  \mathcal{P}(\xi_{t} u^{n-\frac{1}{2}},\psi)  = \displaystyle \frac{1}{\tau \Gamma(2-\alpha)}  \left[ a_{0} \xi_{t} u^{n-\frac{1}{2}} - \displaystyle \sum _{l=1}^{n-1} (a_{n-l-1}-a_{n-l}) \xi_{t} u^{l-\frac{1}{2}} - a_{n-1} \psi  \right] \\
= \displaystyle \frac{a_{0}}{\tau^{2} \Gamma(2-\alpha)} u^{n} - \displaystyle \frac{1}{\tau^2 \Gamma(2-\alpha)} \left[ a_{0} u^{n-1} + \displaystyle \sum _{l=1}^{n-1} (a_{n-l-1}-a_{n-l}) (u^{l}-u^{l-1})  + a_{n-1}\psi  \right].
\end{array}
\end{equation}
Substituting Eq. (\ref{formula 17}) in Eq. (\ref{formula 16}) gives
\begin{equation*}
\begin{array}{l}
\displaystyle \frac{a_{0}}{\tau^{2} \Gamma(2-\alpha)} u^{n} - \displaystyle \frac{1}{\tau^2 \Gamma(2-\alpha)} \left[ a_{0} u^{n-1} + \displaystyle \sum _{l=1}^{n-1} (a_{n-l-1}-a_{n-l}) (u^{l}-u^{l-1})  + a_{n-1}\psi  \right] = \kappa \nabla^{2}u^{n} + g^{n- \frac {1}{2}},
\end{array}
\end{equation*}
by rearranging the above relation we obtain
\begin{equation}\label{formula 18}
\begin{array}{l}
\kappa \nabla^{2}u^{n} + b u^{n}=F^{n-1},
\end{array}
\end{equation}
where
\begin{equation}
b= - \displaystyle \frac{a_{0}}{\tau^{2} \Gamma(2-\alpha)},
\end{equation}
and
\begin{equation}
F^{n-1}=\displaystyle \frac{-1}{\tau^2 \Gamma(2-\alpha)} \left[ a_{0} u^{n-1} + \displaystyle \sum _{l=1}^{n-1} (a_{n-l-1}+a_{n-l}) (u^{l}-u^{l-1}) + a_{n-1}\psi  \right] -
g^{n- \frac {1}{2}}.
\end{equation}
In other words, solving Eq. (\ref{formula 1}) in each time step deduces to solving the above Helmholtz equation. Note that in Eq. (\ref{formula 18}) $b$ is constant and $F^{n-1}$ varies with each time step. In the forthcoming subsection, we employ the constant BEM for solving the Helmholtz equation (\ref{formula 18}).
\section{Spatial discretization via BEM}
In this section we present an overview on how to solve the Helmholtz equation with the BEM. Consider the equation \\
\begin{equation}\label{formula 20}
\begin{array}{l}
\nabla ^{2}u-k^{2} u = \omega, \ \ \ \ \ \ X \in \Omega,\\
u=\bar{u},\ \ \ \ \ \ \ \ \ \ \ \ \ \ \ \ \ X \in \Gamma,
\end{array}
\end{equation}
where $\Omega$ is a bounded computational domain, $\omega$ is a real function and $k$ is a constant. Suppose that $u$ is an unknown function which satisfies Eqs. (\ref{formula 20}). The use of Green's second identity gives
\begin{equation}\label{formula 21}
\begin{array}{l}
\displaystyle \int _{\Omega} \left( \nabla ^{2} u G^{i} - \nabla ^{2} G^{i} u \right) d\Omega = \displaystyle \int _{\Gamma} \left( \displaystyle \frac{\partial u}{\partial n} G^{i} - \frac{\partial G^{i}}{\partial n} u \right) d \Gamma ,
\end{array}
\end{equation}
where $G^{i}=G^{i}(x,y)$ is the fundamental solution of Helmholtz equation based on the source point $P_{i}=(x_i,y_i)$ which satisfies
\begin{equation}\label{formula 22}
\begin{array}{l}
\nabla ^{2}G^{i}-k^{2} G^{i} = \delta (x-P_{i}), \ \ \ \ \ \ X \in \Omega,
\end{array}
\end{equation}
where $\delta(.)$ is the Dirac delta function.
The fundamental solution for two--dimensional Helmholtz equation based on the source point $P_{i}$, is reported as \cite{P. K. Kythe}:
\begin{equation*}
G^{i}(x,y)=\frac{-1}{2 \pi } K_{0}(kr),
\end{equation*}
where $r$ is the distance between the field point $X$ and the source point $P_{i}$ and $K_{0}$ is the modified Bessel function.
Substituting Eqs. (\ref{formula 20}) and (\ref{formula 22}) in Eq. (\ref{formula 21}) we obtain
\begin{equation}\label{formula 23}
\begin{array}{l}
\displaystyle \int_{\Gamma} \left( \displaystyle \frac{\partial u}{\partial n} G^{i} - \displaystyle \frac{\partial G^{i}}{\partial n} u \right) d\Gamma = \displaystyle \int _{\Omega} \left( \nabla^{2} u G^{i} - \nabla^{2} G^{i} u \right) d\Omega \\
\ \ \ \ \ \ \ \ \ \ \  \ \ \ \ \ \ \ \ \ \  \ \ \ \ = \displaystyle \int _{\Omega} \left( (\omega+k^2 u)G^{i} - (\delta(x-P_{i})+k^2G^{i})u \right) d\Omega\\
\ \ \ \ \ \ \ \ \ \ \  \ \ \ \ \ \ \ \ \ \  \ \ \ \ = \displaystyle
\int_{\Omega} \omega G^{i} d\Omega - \displaystyle \frac{\beta}{2\pi}u(P_{i}),
\end{array}
\end{equation}
where $\beta$ is the internal angle regarding to the source point $P_{i}$ for more details see \cite{J. Katsikadelis,Z. Sedaghatjoo}.
To numerical implementation of Eqs. (\ref{formula 1})--(\ref{formula 3}) with BEM, we can write
\begin{equation}\label{formula 24}
\displaystyle \frac{\beta}{2 \pi}u(P_{i})+ \displaystyle \int_{\Gamma} \displaystyle \frac{\partial u}{\partial n}G^{i} d\Gamma -\displaystyle \int_{\Gamma} \displaystyle \frac{\partial G^{i}}{\partial n}u d\Gamma = \displaystyle \int_{\Omega} \omega G^{i} d\Omega, \ \ \ \ p_i \in \Gamma \cup \Omega.
\end{equation}
For approximating the boundary $\Gamma$, we discretize it in to $N$ straight boundary elements $ \displaystyle \{ \Gamma_{i}  \displaystyle \}_{i=1}^{N}$, and we will choose the boundary source point $P_{i}$ at the center of $\Gamma_{i}$ for $i=1,2,...,N$. In constant BEM, the boundary values $u$ and $q= \displaystyle  \frac{\partial u}{\partial n}$ have been approximated constantly over the boundary elements $ \displaystyle \{ \Gamma_{i}  \displaystyle \}_{i=1}^{N}$ equal to $ \displaystyle \{ u_{i}  \displaystyle \}_{i=1}^{N}$ and $ \displaystyle \{ q_{i} \displaystyle \}_{i=1}^{N}$, respectively, where $u_{i}=u(P_{i})$ and $q_{i}=q(P_{i})$ \cite{J. T. Katsikadelis}. So after choosing $L$ interior source points $ \displaystyle \{ P_{i+N}  \displaystyle \}_{i=1}^{L}$ appropriately in $\Omega$, Eq. (\ref{formula 20}) can be approximated as
\begin{equation}\label{formula 25}
\begin{array}{l}
\displaystyle \frac{\beta}{2\pi}u(P_{i}) + \displaystyle \sum_{i=1}^{N} q_{i} \displaystyle \int_{\Gamma_{i}} G^{i} d\Gamma_{i} - \displaystyle \sum_{i=1}^{N} u_{i} \displaystyle \int_{\Gamma_{i}} \displaystyle \frac{\partial G^{i}}{\partial n} d \Gamma_{i} = \displaystyle \int_{\Omega} \omega G^{i} d\Omega, \ \ \ \
i=1,2,...,N.
\end{array}
\end{equation}
Let
\begin{equation}\label{formula 26}
\begin{array}{l}
\bar{H}_{ik}=\displaystyle \int_{\Gamma_{k}} \displaystyle \frac{\partial G^{i}}{\partial n} ds \ \ \ \ , \ \ \ \ G_{ik}=\displaystyle \int_{\Gamma_{k}} G^{i} ds, \ \ \ \ i=1,2,...,N.
\end{array}
\end{equation}
Since we have put the source points (coincide with the mid points of the boundary elements) in Eq. (\ref{formula 25}) as source points, we get the following matrix equation
\begin{equation}\label{formula 27}
\begin{array}{l}
\bf {H} \bf{u^{T}} +  \bf{F} = \bf{G} \bf{q^{T}},
\end{array}
\end{equation}
where $ \bf {H}_{ij}= \bar{H}_{ij} - \displaystyle \frac{1}{2} \delta_{ij}$ ($\delta_{ij}$ is the Kronecker delta function which is defined as $\delta_{ij}=0$ for $i \neq j$ and
$\delta_{ij}=1$ for $i = j$) and $\bf{F}$$= [F^{1},F^{2},...,F^{N}]^{T}$ such that $F^{i}$ has the representation
\begin{equation}\label{formula 28}
\begin{array}{l}
F^{i} =\displaystyle \int_{\Omega} \omega G^{i} d \Omega,  \ \ \ \ i=1,2,...,N.
\end{array}
\end{equation}
In additional $\bf{u}$ and $\bf{q}$ are vector of dimension $N \times 1$.
Now by separating the unknown quantities from the known one, we can obtain the following linear system of equations
\begin{equation}\label{formula 29}
\begin{array}{l}
\bf {AX^{T}} = b^{T}.
\end{array}
\end{equation}
Now we can determine the values of $\bf{u}$ and $\bf{q}$ on the boundary. Therefore after knowing all the boundary quantities on $\Gamma$, the solution $\bf{u}$ can be computed at each internal points $\{ P_{i+N} \}_{i=1}^{L} \in \Omega$ by applying Eq. (\ref{formula 25}) for $\beta=2\pi$.
Hence the solution of Eqs. (\ref{formula 20}) at any interior point $p_i$ can be obtain as follow:
\begin{equation}\label{formula 30}
\begin{array}{l}
u(P_{i}) = \displaystyle \sum_{i=1}^{L} u_{i} G_{i} - \displaystyle \sum_{i=1}^{L} q_{i} H_{i} + \displaystyle \int_{\Omega} \omega G^{i} d\Omega, \ \ \ \
i=1,2,...,L,
\end{array}
\end{equation}
where $\bf{H}$ and $\bf{G}$ can be computed just the same as Eq (\ref{formula 26}) based on the source points $\{ P_{i+N} \}_{i=1}^{L}$.

\subsection{Computation of domain integral}
Suppose $w(X,P)$ is a complex smooth function such that $\nabla^{2}w+bw=1$. For example it may be defined as
\begin{equation*}
w(x,P)=  \left\{ \begin{array}{l}
 \displaystyle \frac{1}{b}(1-\cos(\sqrt{b}(x_{1}-P_{1}))\ \ , \ \ \ \ \ \ if \ \ b\neq 0 \ \ , \\
\\
 \displaystyle \frac{1}{2}(x_{1}-P_{1})^{2}\ \ , \ \ \ \  \ \ \ \ \ \ \ \ \ \ \ \ \ \ \ \ o.w\ \ ,
\end{array} \right.
\end{equation*}
when index $1$ denotes the first component of vectors. In this situation we have
\begin{equation}\label{formula 31}
\begin{array}{l}
\displaystyle \int_{\Omega} \omega(X)G d\Omega =  \displaystyle \int_{\Omega} (\omega(X)-\omega(P))G d\Omega + \omega(P) \displaystyle \int_{\Omega}G d\Omega \\
\ \ \ \ \ \ \ \ \ \ \ \ \ \ \ \ \ =  \acute{I}_{\Omega} + \omega(P)  \displaystyle \int_{\Omega} (\nabla^{2}w+bw)G d\Omega \\
\ \ \ \ \ \ \ \ \ \ \ \ \ \ \ \ =   \acute{I}_{\Omega}   +\omega(P) \left(   \displaystyle \frac{\beta}{2\pi}w(P,P)+ \displaystyle \int_{\Gamma} \frac{\partial w}{\partial n}G -  \displaystyle \frac{\partial G}{\partial n} w d\Gamma
 \right) ,
\end{array}
\end{equation}
where
\begin{equation}\label{formula 32}
\begin{array}{l}
\acute{I}_{\Omega}=\displaystyle \int_{\Omega} \bold{\omega}(X,P) d\Omega \ \ \ \ , \ \ \ \
\bold{\omega}(X,P) = (\omega(X)-\omega(P))G(X,P).
\end{array}
\end{equation}
For the two--dimensional problems $\lim _{X\rightarrow P} \bold{\omega}(X,P)=0$, so if we define $\bold{\omega}(P,P)=0$, the function $\bold{\omega}$ will be continuous at $X=P$. So $\acute{I}_{\Omega}$ can be calculated with no singularity.
We can discretize domain $\Omega$ to smaller sub-domains $\Omega_{j}$ and approximate $\bold{\omega}(X,P)$ over them constantly as $\bold{\omega}(P_{j},P)$ when $P_{j}$ is a collocation point located on middle of $\Omega_{j}$. So we can write
\begin{equation}\label{formula 33}
\begin{array}{l}
\acute{I}_{\Omega} \simeq \displaystyle \sum_{j=1}^{L} \bold{\omega}(P_{j},P)S_{j} \ \ \ \ \ , \ \ \ \ S_{j} = \displaystyle \int_{\Omega_{j}} d \Omega_{j}.
\end{array}
\end{equation}
It can be found from Eq. (\ref{formula 33}) that $S_{j}$ is volume of $\Omega_{j}$, but we can interpret it as the weight of collocation point $P_{j}$ in the domain integral approximation (\ref{formula 33}). For more details, we refer the interested readers to \cite{H. Hosseinzadeh1}.

\section{Implementation of the DRBEM}
In this section we apply DRBEM for the solution of Eqs (\ref{formula 1})--(\ref{formula 3}).
We choose $N$ and $L$ source points into domain and on its boundary, respectively. In this step we make an interpolation for approximating the time derivatives and the inhomogeneous terms in the following:
\begin{equation}\label{formula 50}
\begin{array}{l}
\displaystyle \frac{\partial ^{\alpha} u(X,t)}{\partial t^{\alpha}}  - g(X,t) =  \nu(X,t)
= \displaystyle \sum_{i=1}^{N+L} \eta_{i}(t)\mu_{i}(X),
\end{array}
\end{equation}
where $\mu_i = 1+r_i$ is a linear RBF and  $r_i= \displaystyle \sqrt{(x-x_i)^2+(y-y_i)^2}$ for $i=1,...,N+L$.
From (\ref{formula 1}) and (\ref{formula 50}), we  have
\begin{equation}\label{formula 51}
\begin{array}{l}
\nabla ^2 u= \displaystyle \sum_{i=1}^{N+L} \eta_{i}(t)\mu_{i}(x,y).
\end{array}
\end{equation}
Note that we represent the value of the function $\mu_i$ at source point $p_i$ by $\mu_{ji}$ for $i=1,2,...,N+L$, and set $ \Phi$ a $(N+L)\times (N+L)$ matrix with component $\Phi(j,i)=\mu_{ji}$. From Eq. (\ref{formula 50}) we can write
\begin{equation}\label{formula 52}
\begin{array}{l}
\eta_{i}(t)= \displaystyle \sum_{m=1}^{N+L} \Phi_{im}^{-1} \nu_{m}(t),
\end{array}
\end{equation}
where $\nu_{m}(t)=\nu(x_m,y_m,t)$.
The essential feather in DRBEM is to express $\mu_{i}$, which is a function of $r_i$, as a Laplacian of another function $f_i$. Thus $f_i$ is chosen as the solution to:
\begin{equation}\label{formula 53}
\begin{array}{l}
\nabla ^2 f_i =\mu_i.
\end{array}
\end{equation}
The particular solution $f_i$, for linear RBFs is easily determined as  \cite{Dehghan_1}:
\begin{equation}\label{formula 54}
\begin{array}{l}
f_i = \displaystyle \frac{r_{i}^{2}}{4} + \displaystyle \frac{r_{i}^{3}}{9}.
\end{array}
\end{equation}
In this section, we use a special weight function, that is the fundamental solution for the two--dimensional Laplacian operator, defined by the equation:
\begin{equation*}
\begin{array}{l}
\nabla ^2 G^i = \delta (x-x_i , y-y_i),
\end{array}
\end{equation*}
$\delta$ is the Dirac delta function. The fundamental solution for two--dimensional Laplacian operator based on the source point $P_i$
, is reported as:
\begin{equation*}
\begin{array}{l}
G^i(X) = \displaystyle \frac{-1}{2 \pi} \ln(r),
\end{array}
\end{equation*}
where $r$ is the distance between the field point $(x, y)$ and the source point $P_i$.
Substituting (\ref{formula 53}) into (\ref{formula 51}), writing the right residual formulation of Eq. (\ref{formula 51}) and using the second Green's theorem, lead to:
\begin{equation}\label{formula 55}
\begin{array}{l}
\kappa \left( \displaystyle{\beta}{2\pi}  u_i + \displaystyle \int_{\Gamma} \left[  \frac{\partial G^i}{\partial n} u - G^i \frac{\partial u}{\partial n}  \right] d\Gamma  \right)
= \displaystyle \sum_{j=1}^{N+L} \eta_j(t) \left\{ \displaystyle{\beta}{2\pi} f_{ij} +  \int_{\Gamma} \left[  \frac{\partial G^i}{\partial n} f_j -   G^i \frac{\partial f_j}{\partial n}  \right] d\Gamma   \right\} ,
\end{array}
\end{equation}
for $i=1,2,...,N+L$, where $f_{ij}=f_i(x_j)$.
The term $\displaystyle \frac{\partial f_i}{\partial n}$ is the normal derivative of $f_i$ and can be expressed as \cite{Dehghan_1}:
\begin{equation}\label{formula 56}
\begin{array}{l}
\displaystyle \frac{\partial f_i}{\partial n}= \left(  \displaystyle \frac{1}{2} + \displaystyle \frac{r_i}{3}  \right) \left[ (x-x_i)n_x+(y-y_i)n_y \right],
\end{array}
\end{equation}
where $n=(n_x,n_y)$ is the unit outward normal to the boundary $\Gamma$. Next we show the discretized form of Eq. (\ref{formula 55}), with summations over the constant boundary elements, for the $i$ th source point as follow
\begin{equation}\label{formula 57}
\begin{array}{l}
\kappa \left\{ \displaystyle \frac{\beta}{2\pi} u_i + \displaystyle \sum_{q=1}^{N} \displaystyle \int_{\Gamma_q} \left[  \frac{\partial G^i}{\partial n} u - G^i \frac{\partial u}{\partial n}  \right] d\Gamma \right\}
= \displaystyle \sum_{j=1}^{N+L} \eta_j(t) \left( \displaystyle \frac{\beta}{2\pi} f_{ij} +  \displaystyle \sum_{q=1}^{N} \int_{\Gamma_q} \left[  \frac{\partial G^i}{\partial n} f_j -   G^i \frac{\partial f_j}{\partial n}  \right] d\Gamma   \right) .
\end{array}
\end{equation}
The expression given by (\ref{formula 57}) can be expressed in terms of the coefficients $\bar{ H }_{iq}$ and $G_{iq}$, using the
constant element approximation \cite{J. Katsikadelis}. The resulting equation is
\begin{equation}\label{formula 58}
\begin{array}{l}
\kappa \left( \displaystyle \frac{\beta}{2\pi} u_i + \displaystyle \sum_{q=1}^{N}  \bar{H}_{iq} u_q - \displaystyle \sum_{q=1}^{N} G_{iq} \frac{\partial u_q}{\partial n} \right)
= \displaystyle \sum_{j=1}^{N+L}  \eta_j(t) \left\{ \displaystyle \frac{\beta}{2\pi} f_{ij} +
\displaystyle \sum_{q=1}^{N} \bar{H}_{iq} f_{qj} -  \displaystyle \sum_{q=1}^{N} G_{iq} \frac{\partial f_{qj}}{\partial n}\right\} ,
\end{array}
\end{equation}
for $i=1,2,...,N$ (we defined $\bar{ H }$ and $G$ in Eq. (\ref{formula 26})).
For simplicity we put, $\bold{q}=\displaystyle \frac{\partial u}{\partial n}$, $\hat{\bold{u}} = f$ and $\hat{\bold{q}} = \displaystyle \frac{\partial f}{\partial n}$.
Eq. (\ref{formula 58}) can be expressed in the matrix form as
\begin{equation}\label{formula 60}
\begin{array}{l}
\kappa \left(  \bold{H} \bold{u} - \bold{G} \bold{q} \right)
= \displaystyle \sum_{j=1}^{N+L} \eta_j(t) \left\{ \bold{H} \hat{\bold{u}}_j - \bold{G} \hat{\bold{q}}_j \right\}.
\end{array}
\end{equation}
If each of the vectors $\hat{\bold{u}}_k$ and $\hat{\bold{q}}_k$ are considered to be one column of the matrices $\hat{\bold{U}}$ and $\hat{\bold{Q}}$, respectively, then for the solution of the boundary values we have

\begin{equation}\label{formula 61}
\begin{array}{l}
\kappa \left( \bold{H}^{bs} \bold{u} - \bold{G}^{bs} \bold{q} \right)
= \left( \bold{H}^{bs} \hat{\bold{U}}^{bs} - \bold{G}^{bs} \hat{\bold{Q}}^{bs} \right) \bold{ \eta }.
\end{array}
\end{equation}
Eq. (\ref{formula 61}) is the basis of the application of the DRBEM and involves discretization of the boundary only. One can calculate the values of $\bold{u}$ at the internal nodes, by redefining the matrices $\bold{\bar{H}}$ and $\bold{G}$ in Eq. (\ref{formula 58}) in such a way that they contain the point $p_i$ as an interior node. So for the internal nodes we have
\begin{equation}\label{formula 62}
\begin{array}{l}
\kappa \left( \bold{H}^{is} \bold{u} -\bold{I}\bold{u}  - \bold{G}^{is} \bold{q} \right)
= \left( \bold{H}^{is} \hat{\bold{U}}^{is} -\bold{I} \hat{\bold{U}}^{is} - \bold{G}^{is} \hat{\bold{Q}}^{is} \right) \bold{ \eta },
\end{array}
\end{equation}
where $\eta=[ \eta_1,\eta_2,...,\eta_{N+L} ]^T$ and $bs$ and $is$ refer to boundary and internal solutions, respectively.
Eq. (\ref{formula 52}) yields
\begin{equation}\label{formula 63}
\begin{array}{l}
\bold{\eta} = \bold{\Phi}^{-1} \bold{\nu}.
\end{array}
\end{equation}
From Eqs. (\ref{formula 50}), (\ref{formula 63}) we obtain
\begin{equation}\label{formula 64}
\begin{array}{l}
\eta = \bold{\Phi} ^{-1} ( \displaystyle \frac{\partial^\alpha \bold{u}}{\partial t^{\alpha}}  - \bold{g} ).
\end{array}
\end{equation}
Combining Eqs. (\ref{formula 61}), (\ref{formula 62}) and (\ref{formula 64}) together at one step, one can get the following system of equations \cite{C. Bozkaya}:
\begin{equation}\label{formula 65}
\begin{array}{l}
\kappa \left( \left[ {\begin{array}{*{20}{c}} {\bold{H}^{bs}} \ \ \ \ {\bold{0}}\\ {\bold{H}^{is}} \ \ \ \ {\bold{-I}} \end{array}} \right] \left[ {\begin{array}{*{20}{c}} {\bold{u}^{bs}} \\ {\bold{u}^{is}}  \end{array}} \right]
- \left[ {\begin{array}{*{20}{c}} {\bold{G}^{bs}} \ \ \ \ {\bold{0}}\\
{\bold{G}^{is}} \ \ \ \ {\bold{0}} \end{array}} \right]
\left[ {\begin{array}{*{20}{c}} {\bold {q}^{bs}} \\ {0}  \end{array}} \right]\right)
= \bold{D} \left( \displaystyle \frac{\partial^\alpha \bold{u}}{\partial t^{\alpha}}  - \bold{g} \right) ,
\end{array}
\end{equation}
where
\begin{equation*}
\begin{array}{l}
\bold{D} = \left( \left[ {\begin{array}{*{20}{c}} {\bold{H}^{bs}} \ \ \ \ {\bold{0}}\\ {\bold{H}^{is}} \ \ \ \ {\bold{-I}} \end{array}} \right]
\left\{ {\begin{array}{*{20}{c}} {\hat{U}^{bs}} \\ {\hat{U}^{is}}  \end{array}} \right\}
- \left[ {\begin{array}{*{20}{c}} {\bold{G}^{bs}} \ \ \ \ {\bold{0}}\\ {\bold{G}^{is}} \ \ \ \ {\bold{0}} \end{array}} \right]
\left\{ {\begin{array}{*{20}{c}} {\hat{Q}^{bs}} \\ {0}  \end{array}} \right\} \right)  \bold{\Phi} ^{-1}.
\end{array}
\end{equation*}
Now Eq. (\ref{formula 65}) constitutes a system of $N+L$ equations in $N+L$ unknown functions of $t$.
We use the Caputo sense definition (Eqs. (\ref{formula 6})--(\ref{formula 7}))
and make the following approximation:
\begin{equation}\label{formula 68}
\begin{array}{l}
\bold{g} \simeq \displaystyle \frac{1}{2}(\bold{g}^{k}+\bold{g}^{k-1}).
\end{array}
\end{equation}
For $u$ we employ the subsequent approximation:
\begin{equation}\label{formula 69}
\begin{array}{l}
u(X,t) \simeq \displaystyle \frac{1}{2}(\bold{u}^{k}+\bold{u}^{k-1}).
\end{array}
\end{equation}
For simplicity, let
\begin{equation*}
\bold{\tilde{H}}= \left[ {\begin{array}{*{20}{c}} {\bold{H}^{bs}} \ \ \ \ {\bold{0}}\\ {\bold{H}^{is}} \ \ \ \ {\bold{-I}} \end{array}} \right], \ \ \ \ \bold{\tilde{G}}= \left[ {\begin{array}{*{20}{c}} {\bold{G}^{bs}} \ \ \ \ {\bold{0}}\\ {\bold{G}^{is}} \ \ \ \ {\bold{0}} \end{array}} \right], \ \ \ \
\vartheta = \displaystyle \frac{ 1 }{\tau^2 \Gamma(2-\alpha)},
\end{equation*}
where $\bold{I}$ is $(N+L)\times (N+L)$ identity matrix.
Eqs. (\ref{formula 65})--(\ref{formula 69}) yield :
\begin{equation}\label{formula 70}
\begin{array}{l}
\displaystyle \frac {1}{2} \left\{  \bold{\tilde{H}} \left( \left[ {\begin{array}{*{20}{c}} {\bold{u}^{bs}} \\ {\bold{u}^{is}}  \end{array}} \right]^{k}
+ \left[ {\begin{array}{*{20}{c}} {\bold{u}^{bs}} \\ {\bold{u}^{is}}  \end{array}} \right]^{k-1} \right)
- \bold{\tilde{G}}
\left( \left[ {\begin{array}{*{20}{c}} {\bold{q}^{bs}} \\ {0}  \end{array}} \right]^{k} +\left[ {\begin{array}{*{20}{c}} {\bold{q}^{bs}} \\ {0}  \end{array}} \right]^{k-1} \right) \right\}\\
= \bold{D} \left(
\displaystyle \frac{1}{\tau \Gamma(2-\alpha)} \left[ a_{0} \xi_t \bold{u}^{k-\frac{1}{2}}  - \displaystyle \sum_{l=1}^{k-1} (a_{n-k-1}-a_{n-k})\xi_t \bold{u}^{l-\frac{1}{2}}
- a_{k-1} \psi \right] \right)
- \displaystyle \frac {1}{2} \bold{D} \left( \left[ {\begin{array}{*{20}{c}} {\bold{g}^{bs}} \\ {\bold{g}^{is}}  \end{array}} \right]^{k}
+ \left[ {\begin{array}{*{20}{c}} {\bold{g}^{bs}} \\ {\bold{g}^{is}}  \end{array}} \right]^{k-1} \right) .
\end{array}
\end{equation}
Rearranging Eq. (\ref{formula 70}), the following form can be obtained

\begin{equation}\label{formula 71}
\begin{array}{l}
\left\{  \displaystyle \frac {1}{2} \bold{\tilde{H}} - a_0  \vartheta \bold{D}  \right\} \left[ {\begin{array}{*{20}{c}} {\bold{u}^{bs}} \\ {\bold{u}^{is}}  \end{array}} \right]^{k}
- \displaystyle \frac {1}{2} \bold{\tilde{G}}
\left[ {\begin{array}{*{20}{c}} {\bold{q}^{bs}} \\ {0}  \end{array}} \right]^{k}
= -  \left\{ \displaystyle \frac {1}{2} \bold{\tilde{H}} - a_0  \vartheta \bold{D}  \right\} \left[ {\begin{array}{*{20}{c}} {\bold{u}^{bs}} \\ {\bold{u}^{is}}  \end{array}} \right]^{k-1}
+\displaystyle \frac {1}{2} \bold{\tilde{G}}
\left[ {\begin{array}{*{20}{c}} {\bold{q}^{bs}} \\ {0}  \end{array}} \right]^{k-1}\\
- \bold{D} \left( \vartheta \left[ \displaystyle \sum_{l=1}^{k-1} (a_{k-j-1}-a_{k-j}) \bold{u}^{l-\frac{1}{2}} + a_{k-1} \psi \right] \right) - \displaystyle \frac {1}{2} \bold{D} \left( \left[ {\begin{array}{*{20}{c}} {\bold{g}^{bs}} \\ {\bold{g}^{is}}  \end{array}} \right]^{k}
+ \left[ {\begin{array}{*{20}{c}} {\bold{g}^{bs}} \\ {\bold{g}^{is}}  \end{array}} \right]^{k-1} \right) .
\end{array}
\end{equation}
It is necessary to separate the unknown from the known quantities. After moving all the known quantities to the left hand side, we obtain the following linear system of equations
\begin{equation}\label{formula 29}
\begin{array}{l}
\bold{A}\bold{x} = \bold{c},
\end{array}
\end{equation}
where the vector $\bold{x}$ contains $N$ unknown boundary values $\bold{q}$ and $L$ unknown internal values of $\bold{u}$.

\section{Numerical results}
In this section we present the numerical results of the proposed methods on several test problems. We test the accuracy and the stability
of the method described in this paper by performing the mentioned methods for different values of $N$, $\tau$ and $\alpha$.
To see the convergence of the proposed method we employ the following error norm
\begin{equation*}
RMS--error =\sqrt{ \frac{1}{L} \displaystyle \sum_{i=1}^{L} ( u_i-\tilde{u}_i )^2},
\end{equation*}
where $u$ and $\tilde{u}$ denote the exact and approximated solution, respectively and $L$ denotes the number of internal points.

\subsection{Test problem 1}
Consider the following problem \cite{Ya-Nan Zhang}:
\begin{equation*}
\begin{array}{l}
D_t^{\alpha}u(x,y,t)= \nabla^{2} u(x,y,t) + f(x,y,t), \ \ (x,y) \in \Omega, \ \ 0 < t \leq 1,
\end{array}
\end{equation*}
with boundary condition
\begin{equation*}
\begin{array}{l}
u(x,y,t)=0, \ \ (x,y) \in \partial \Omega, \ \ \ \ 0 < t \leq 1,
\end{array}
\end{equation*}
and initial condition
\begin{equation*}
\begin{array}{l}
u(x,y,0) = 0, \ \ \psi(x,y)=0, \ \ \ \ (x,y) \in \Omega,
\end{array}
\end{equation*}
where $\Omega=[0,\pi]^{2}$. Then we choose the exact solution:
\begin{equation*}
\begin{array}{l}
u(x,y,t)= t^{2+\alpha} \sin( x) \sin (y).
\end{array}
\end{equation*}
It can be checked that the associated forcing term is
\begin{equation*}
\begin{array}{l}
f(x,y,t)= \sin (x) \sin (y) \left[ \displaystyle \frac{\Gamma(3+\alpha)}{2} t^{2} +2 t^{2+\alpha} \right].
\end{array}
\end{equation*}
We solve this problem with the methods presented in this article with
several values of $\tau$, $N$, $\alpha$ for final time $T=1$. The estimated errors applied methods are shown in Table 1. In Table 2 the accuracy of the methods is shown for several values of $N$.
Fig. 1 shows the absolute error with $N=100$, $\tau=1/100$ and $\alpha=1.05$ at final time $T=1$.

\begin{center}
\begin{tabular}{llllllllllll}
\multicolumn{12}{c}{\footnotesize{\textbf{Table 1}}}\\
\multicolumn{12}{c}{\footnotesize{\textbf{Numerical stability with $N = 80$ for Test problem 1}}}\\
\hline
\hline
\cline{4-6}\cline{8-10}
&\footnotesize{$\alpha$}&&\footnotesize{$\tau$}&&\footnotesize{$BEM$}&&\footnotesize{}&&\footnotesize{$DRBEM$}&&\\
\hline
\hline
&\footnotesize{$1.25$}&&\footnotesize{$1/2$}&&\footnotesize{$1.4937\times10^{-2}$}&&\footnotesize{}&&\footnotesize{$7.3113\times10^{-2}$}&&\\
&\footnotesize{}&&\footnotesize{$1/4$}&&\footnotesize{$7.6730\times10^{-3}$}&&\footnotesize{}&&\footnotesize{$2.2348\times10^{-2}$}&&\\
&\footnotesize{}&&\footnotesize{$1/8$}&&\footnotesize{$3.8253\times10^{-3}$}&&\footnotesize{}&&\footnotesize{$6.7127\times10^{-3}$}&&\\
&\footnotesize{}&&\footnotesize{$1/16$}&&\footnotesize{$1.8655\times10^{-3}$}&&\footnotesize{}&&\footnotesize{$1.9858\times10^{-3}$}&&\\
\hline
\hline
&\footnotesize{$1.75$}&&\footnotesize{$1/2$}&&\footnotesize{$3.4845\times10^{-2}$}&&\footnotesize{}&&\footnotesize{$3.0726\times10^{-1}$}&&\\
&\footnotesize{}&&\footnotesize{$1/4$}&&\footnotesize{$1.7287\times10^{-2}$}&&\footnotesize{}&&\footnotesize{$1.2868\times10^{-1}$}&&\\
&\footnotesize{}&&\footnotesize{$1/8$}&&\footnotesize{$8.2701\times10^{-3}$}&&\footnotesize{}&&\footnotesize{$5.3983\times10^{-2}$}&&\\
&\footnotesize{}&&\footnotesize{$1/16$}&&\footnotesize{$3.8777\times10^{-3}$}&&\footnotesize{}&&\footnotesize{$2.2705\times10^{-2}$}&&\\
\hline
\hline
\end{tabular}
\end{center}

\begin{center}
\begin{tabular}{llllllllllll}
\multicolumn{12}{c}{\footnotesize{\textbf{Table 2}}}\\
\multicolumn{12}{c}{\footnotesize{\textbf{ The results of Test problem 1}}}\\
\hline
\hline
\cline{4-6}\cline{8-10}
&\footnotesize{}&&\footnotesize{}&&\footnotesize{$\tau=1/20$}&&\footnotesize{$\tau=1/500$}&&\footnotesize{}&&\footnotesize{}\\
\cline{5-6}\cline{8-10}
&\footnotesize{$\alpha$}&&\footnotesize{$N$}&&\footnotesize{$BEM$}&&\footnotesize{$DRBEM$}&&\\
\hline
\hline
&\footnotesize{$1.25$}&&\footnotesize{$20$}&&\footnotesize{$8.6823\times10^{-4}$}&&\footnotesize{$2.2966\times10^{-3}$}&&\\
&\footnotesize{}&&\footnotesize{$40$}&&\footnotesize{$2.0403\times10^{-4}$}&&\footnotesize{$3.6969\times10^{-4}$}&&\\
&\footnotesize{}&&\footnotesize{$80$}&&\footnotesize{$9.7308\times10^{-5}$}&&\footnotesize{$4.6334\times10^{-5}$}&&\\
&\footnotesize{}&&\footnotesize{$160$}&&\footnotesize{$1.8993\times10^{-5}$}&&\footnotesize{$1.0820\times10^{-5}$}&&\\
\hline
\hline
&\footnotesize{$1.5$}&&\footnotesize{$20$}&&\footnotesize{$1.7767\times10^{-3}$}&&\footnotesize{$2.0063\times10^{-3}$}&&\\
&\footnotesize{}&&\footnotesize{$40$}&&\footnotesize{$3.3122\times10^{-4}$}&&\footnotesize{$2.3330\times10^{-4}$}&&\\
&\footnotesize{}&&\footnotesize{$80$}&&\footnotesize{$1.5129\times10^{-4}$}&&\footnotesize{$4.1112\times10^{-5}$}&&\\
&\footnotesize{}&&\footnotesize{$160$}&&\footnotesize{$6.7731\times10^{-6}$}&&\footnotesize{$4.4468\times10^{-5}$}&&\\
\hline
\hline
\end{tabular}
\end{center}

\subsection{Test problem 2}
Consider the following problem:
\begin{equation*}
D_t^{\alpha}u(x,y,t)= \nabla^{2} u(x,y,t) + 2 \exp(x+y)(t^{2}-\displaystyle \frac {t^{2-\alpha}}{\Gamma(3-\alpha)}),
\end{equation*}
in a circle with center $(0,0)$ and radius 1, with Dirichlet boundary condition, and initial boundary condition $u^{0}(x,y)=0$ (Fig. 4 (right plane) shows the computational domain).
The exact solution of this test problem which is taken from \cite{A. Mohebbi1} is:
\begin{equation*}
u(x,y,t)=\exp(x+y)t^{2}.
\end{equation*}

\begin{center}
\begin{tabular}{llllllllllll}
\multicolumn{12}{c}{\footnotesize{\textbf{Table 3}}}\\
\multicolumn{12}{c}{\footnotesize{\textbf{ The results of Test problem 2 with $\tau= 1/20$ \ \ \ \ \ \ \ \ \ \ }}}\\
\hline
\hline
&\footnotesize{$\alpha$}&&\footnotesize{$N$}&&\footnotesize{$BEM$}&&\footnotesize{$DRBEM$}&&\\
\hline
\hline
&\footnotesize{$1.25$}&&\footnotesize{$20$}&&\footnotesize{$3.9578\times10^{-3}$}&&\footnotesize{$9.6255\times10^{-3}$}&&\\
&\footnotesize{}&&\footnotesize{$40$}&&\footnotesize{$8.0834\times10^{-4}$}&&\footnotesize{$2.2821\times10^{-3}$}&&\\
&\footnotesize{}&&\footnotesize{$80$}&&\footnotesize{$3.1444\times10^{-4}$}&&\footnotesize{$4.9215\times10^{-4}$}&&\\
&\footnotesize{}&&\footnotesize{$160$}&&\footnotesize{$1.4864\times10^{-4}$}&&\footnotesize{$1.2339\times10^{-4}$}&&\\
\hline
\hline
&\footnotesize{$1.75$}&&\footnotesize{$20$}&&\footnotesize{$1.5021\times10^{-2}$}&&\footnotesize{$9.5778\times10^{-3}$}&&\\
&\footnotesize{}&&\footnotesize{$40$}&&\footnotesize{$1.3229\times10^{-3}$}&&\footnotesize{$2.2722\times10^{-3}$}&&\\
&\footnotesize{}&&\footnotesize{$80$}&&\footnotesize{$3.5703\times10^{-4}$}&&\footnotesize{$4.8886\times10^{-4}$}&&\\
&\footnotesize{}&&\footnotesize{$160$}&&\footnotesize{$1.6264\times10^{-4}$}&&\footnotesize{$1.2227\times10^{-4}$}&&\\
\hline
\hline
\end{tabular}
\end{center}

We solve this problem with several values of $N$, $\tau$ and $\alpha$ for final time $T=1$. The error of applied methods are shown in Table 3.
Fig. 2 shows the absolute error with $\alpha=1.95$ for some different values of final time $T$. According this figure, as final time $T$ increases, the accuracy of DRBEM is satisfactory.
\subsection{Test problem 3}
We consider the linear TFPDE:
\begin{equation}\label{formula 34}
\begin{array}{l}
D_t^{\alpha}u(x,y,t)= \nabla^{2} u(x,y,t) -2\cos(\pi(x+y))(\pi^{2}t^{2}+\displaystyle \frac{ t^{2-\alpha }}{\Gamma(3-\alpha)}), \ \ \ \ (x,y) \in \Omega, \ \ \ \ (x,y) \in \Omega,
\end{array}
\end{equation}
with the initial condition
\begin{equation*}
\begin{array}{l}
u(x,y,0)=0, \ \ \ \ \ \ (x,y) \in \Omega ,
\end{array}
\end{equation*}
and the Dirichlet boundary conditions. The exact solution of the problem is
\begin{equation*}
\begin{array}{l}
u(x,y,t)=\cos(\pi (x+y)) t^{2},
\end{array}
\end{equation*}
where $\Omega$ is a bounded computational domain shown in Fig. 4 (left plane).
The error of applied methods are shown in Table 4. Fig. 3 shows the absolute error with $\tau=1/10$ for several values of $\alpha$ and $N$ at final time $T=1$.
\begin{center}
\begin{tabular}{llllllllllll}
\multicolumn{12}{c}{\footnotesize{\textbf{Table 4}}}\\
\multicolumn{12}{c}{\footnotesize{\textbf{ The results of Test problem 3 with $\tau=1/10$}}}\\
\hline
\hline
\cline{4-6}\cline{8-10}
&\footnotesize{$\alpha$}&&\footnotesize{$N$}&&\footnotesize{$BEM$}&&\footnotesize{}&&\footnotesize{$DRBEM$}&&\\
\hline
\hline
&\footnotesize{$1.25$}&&\footnotesize{$50$}&&\footnotesize{$7.4715\times10^{-2}$}&&\footnotesize{}&&\footnotesize{$5.5703\times10^{-2}$}&&\\
&\footnotesize{}&&\footnotesize{$100$}&&\footnotesize{$3.1354\times10^{-2}$}&&\footnotesize{}&&\footnotesize{$2.0639\times10^{-2}$}&&\\
&\footnotesize{}&&\footnotesize{$200$}&&\footnotesize{$1.5387\times10^{-2}$}&&\footnotesize{}&&\footnotesize{$7.7330\times10^{-3}$}&&\\
&\footnotesize{}&&\footnotesize{$400$}&&\footnotesize{$7.6845\times10^{-3}$}&&\footnotesize{}&&\footnotesize{$2.6892\times10^{-3}$}&&\\
\hline
\hline
&\footnotesize{$1.75$}&&\footnotesize{$50$}&&\footnotesize{$7.0550\times10^{-2}$}&&\footnotesize{}&&\footnotesize{$5.1564\times10^{-2}$}&&\\
&\footnotesize{}&&\footnotesize{$100$}&&\footnotesize{$2.2724\times10^{-2}$}&&\footnotesize{}&&\footnotesize{$1.9106\times10^{-2}$}&&\\
&\footnotesize{}&&\footnotesize{$200$}&&\footnotesize{$1.0464\times10^{-2}$}&&\footnotesize{}&&\footnotesize{$7.1363\times10^{-3}$}&&\\
&\footnotesize{}&&\footnotesize{$400$}&&\footnotesize{$5.2017\times10^{-3}$}&&\footnotesize{}&&\footnotesize{$2.4746\times10^{-3}$}&&\\
\hline
\hline
\end{tabular}
\end{center}

\section{Conclusion}
This paper has proposed a boundary elements method (BEM) and dual reciprocity boundary elements method (DRBEM) for the numerical solution of two--dimensional time-fractional diffusion--wave. We have shown that the approximation of solution in each time step has been reduced to solving an inhomogeneous Helmholtz equation in BEM. The inhomogeneous part of the mentioned Helmholtz equation in each time steps, depends on all previous time steps. This inhomogeneous Helmholtz equation is approximated with the constant BEM. The main concerning problem of the BEM is the domain integrals that appear in the integral representation of this Helmholtz equation. So, in each time step we need to calculate a domain integral. Since these domain integrals are singular, the singularity of these domain integrals can be omitted via the simple technique. The mentioned technique has been also presented in this paper. These domain integrals have been approximated by Monte--Carlo approach. Also the time--fractional derivative of problems is described in the Caputo sense.
The DRBEM is applied to eliminate the domain integrals which appear in the boundary integral equation. The linear radial basis function (RBF) was used for interpolation of the nonlinear, inhomogeneous and time derivative terms. Some numerical tests taken from the literature are considered on regular and complex regions which confirm the validity of the presented techniques.


\newpage
\begin{figure}[t]
\includegraphics[width=8cm, height=7.cm]{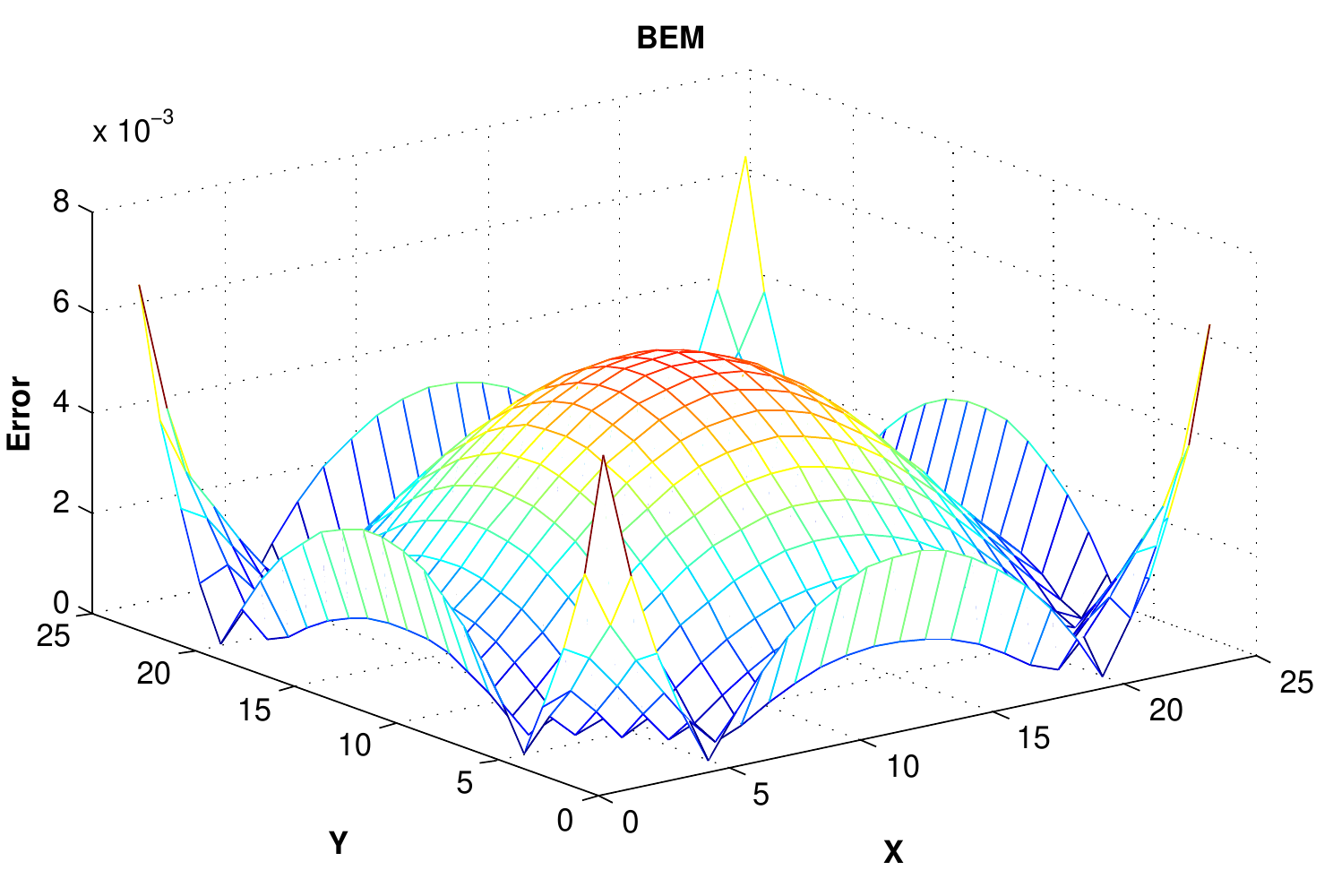}
\includegraphics[width=8cm, height=7.cm]{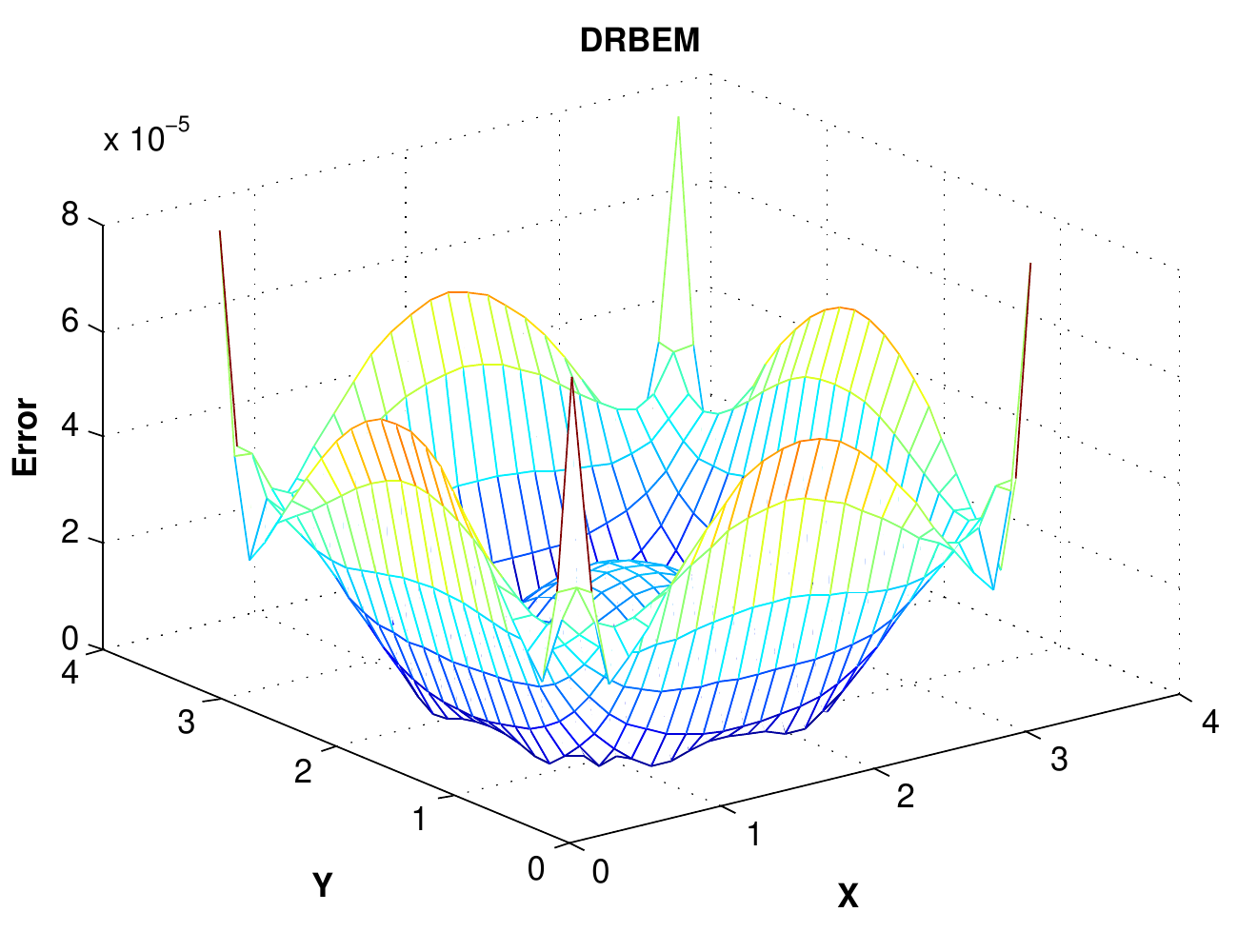}
\caption{Graphs of the average values of relative error with $\tau=1/100$, $N=100$ and $\alpha=1.05$ at final time $T=1$ for Test problem 1.}
\end{figure}

\begin{figure}[t]
\includegraphics[width=8cm, height=7.cm]{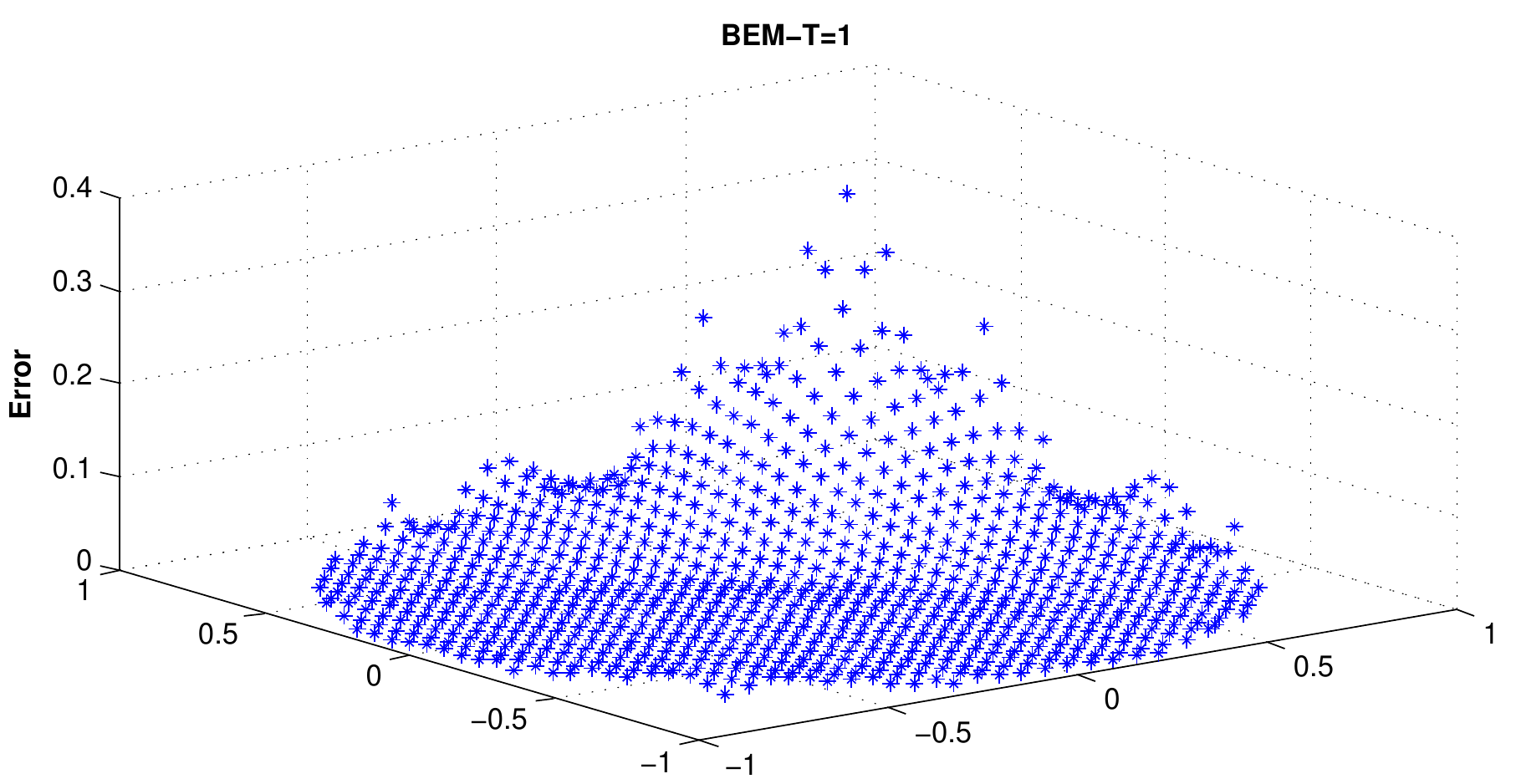}
\includegraphics[width=8cm, height=7.cm]{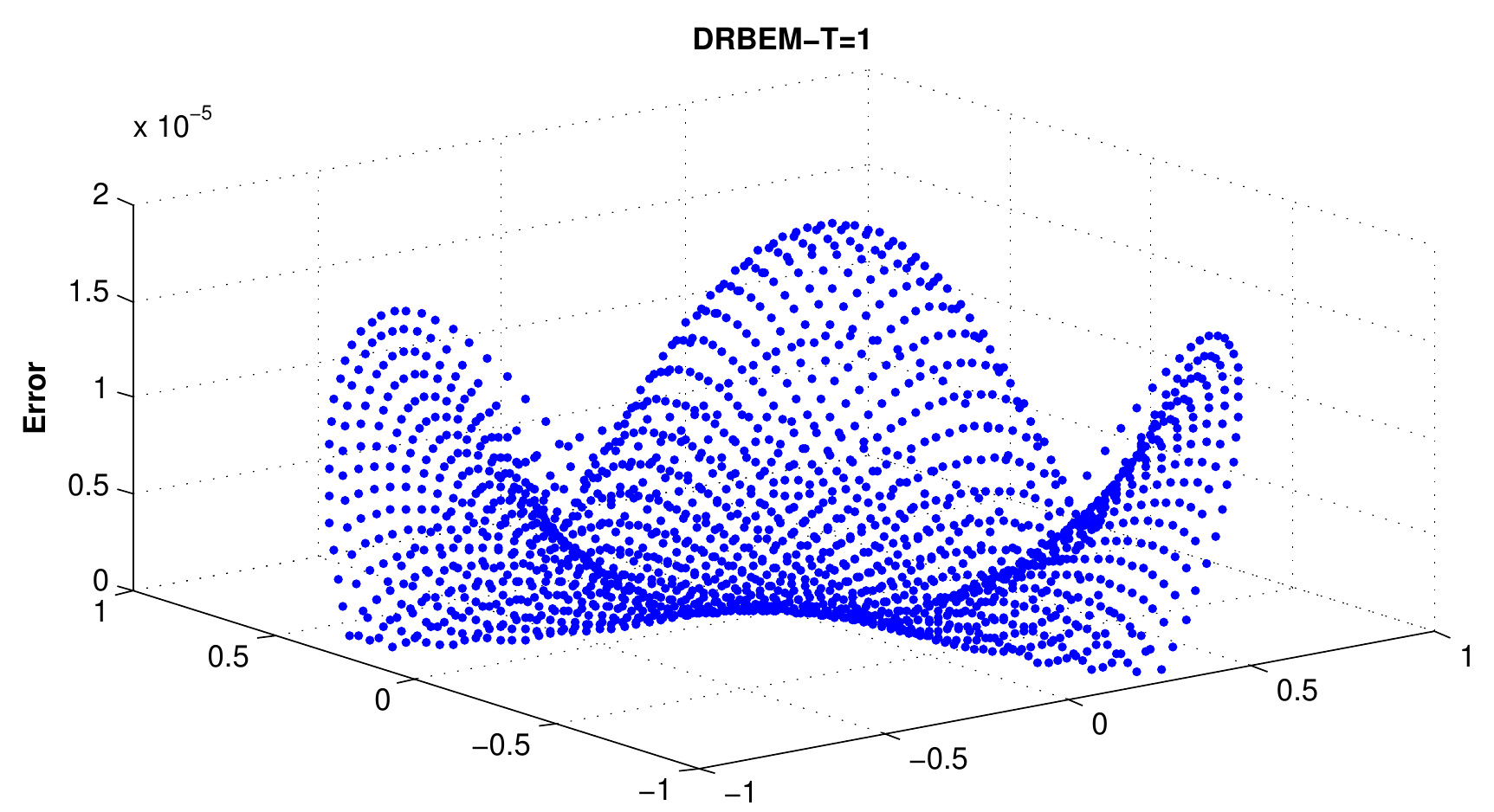}
\includegraphics[width=8cm, height=7.cm]{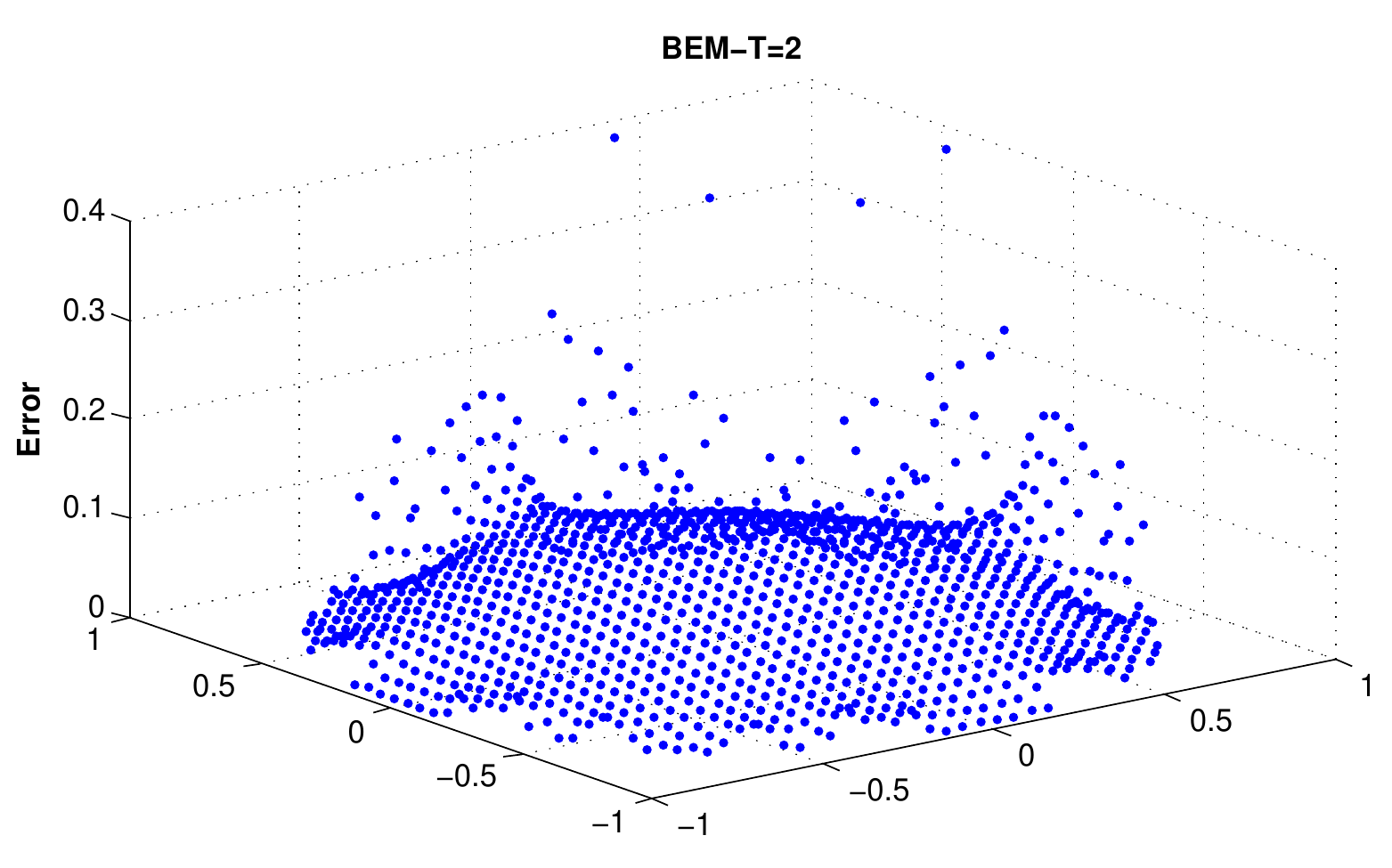}
\includegraphics[width=8cm, height=7.cm]{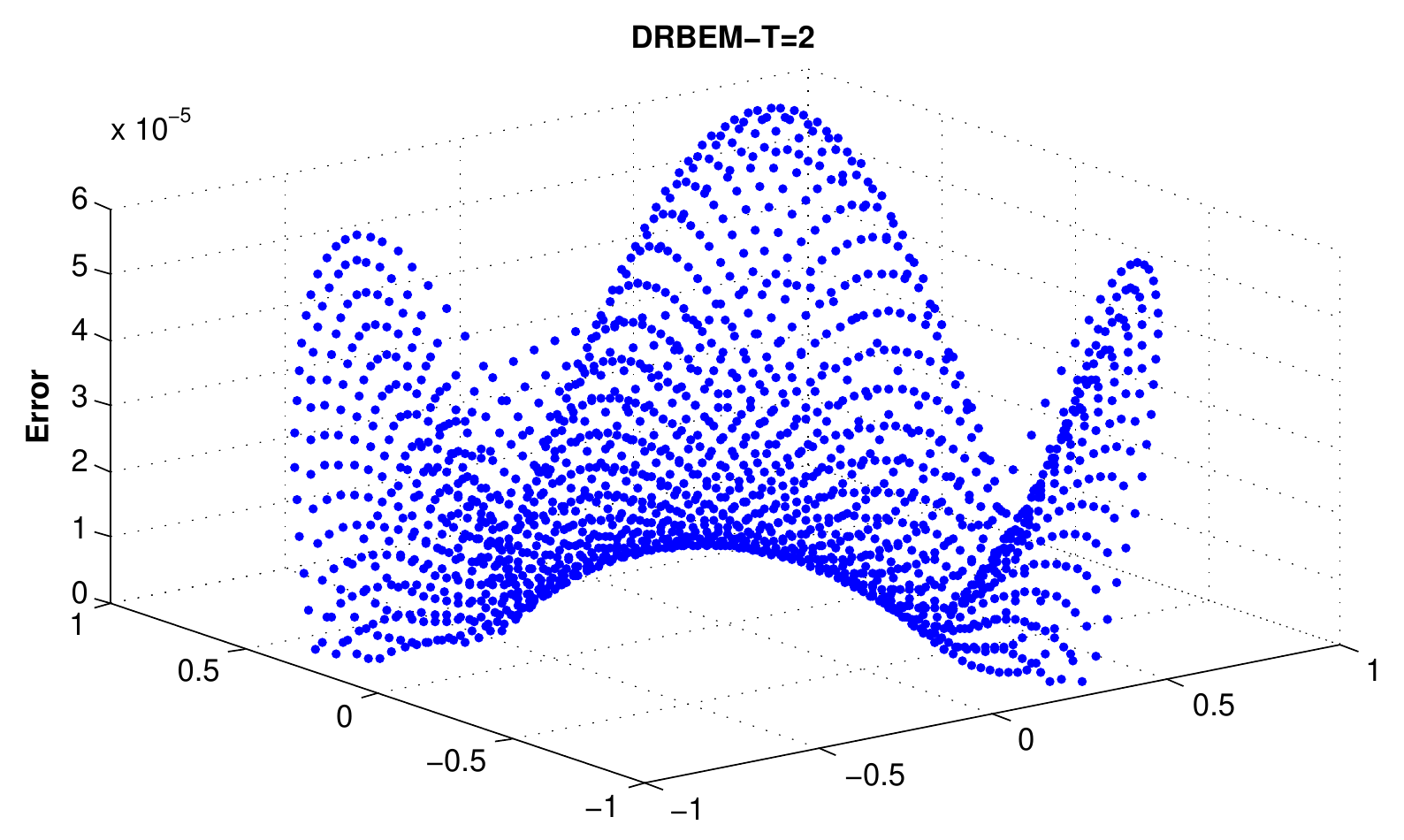}
\includegraphics[width=8cm, height=7.cm]{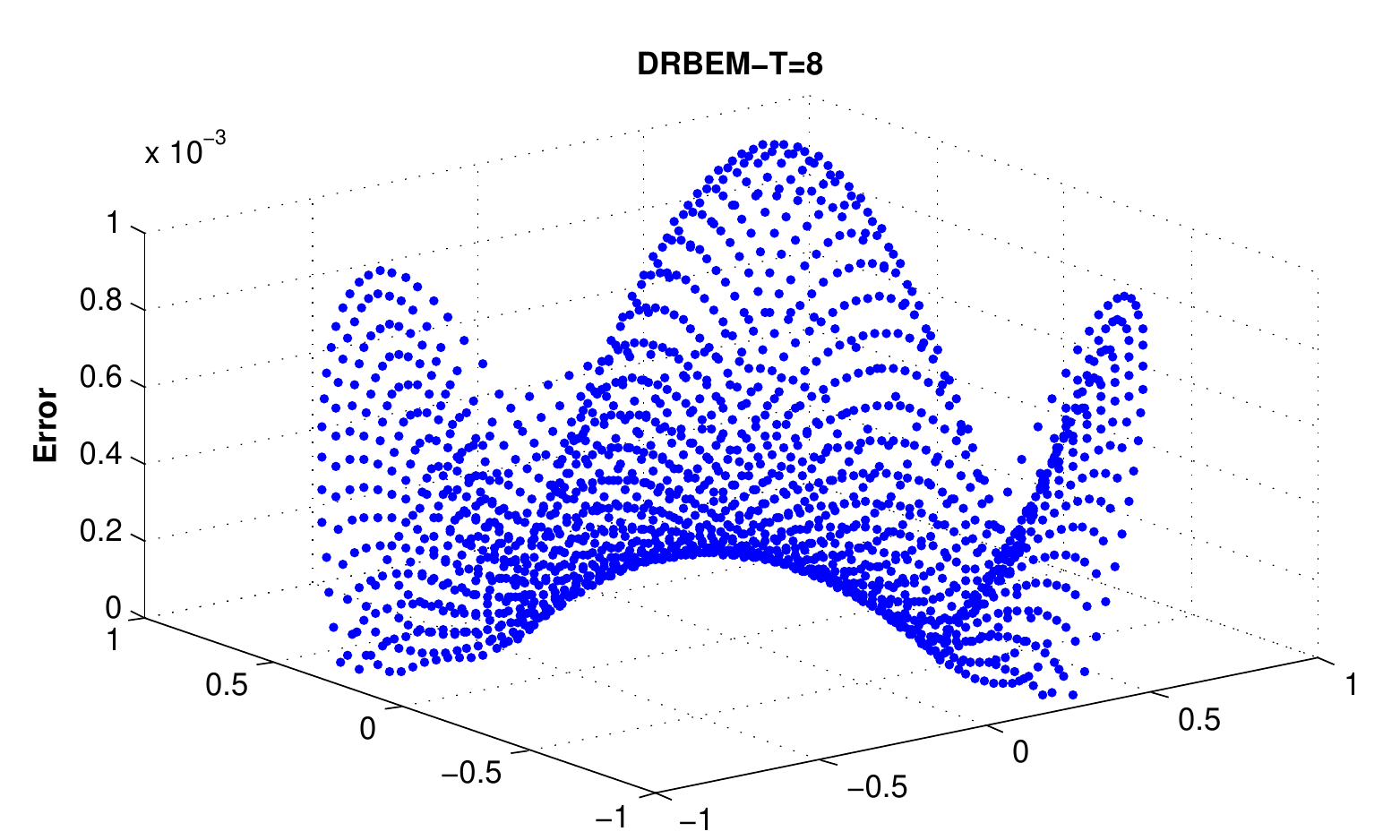}
\includegraphics[width=8cm, height=7.cm]{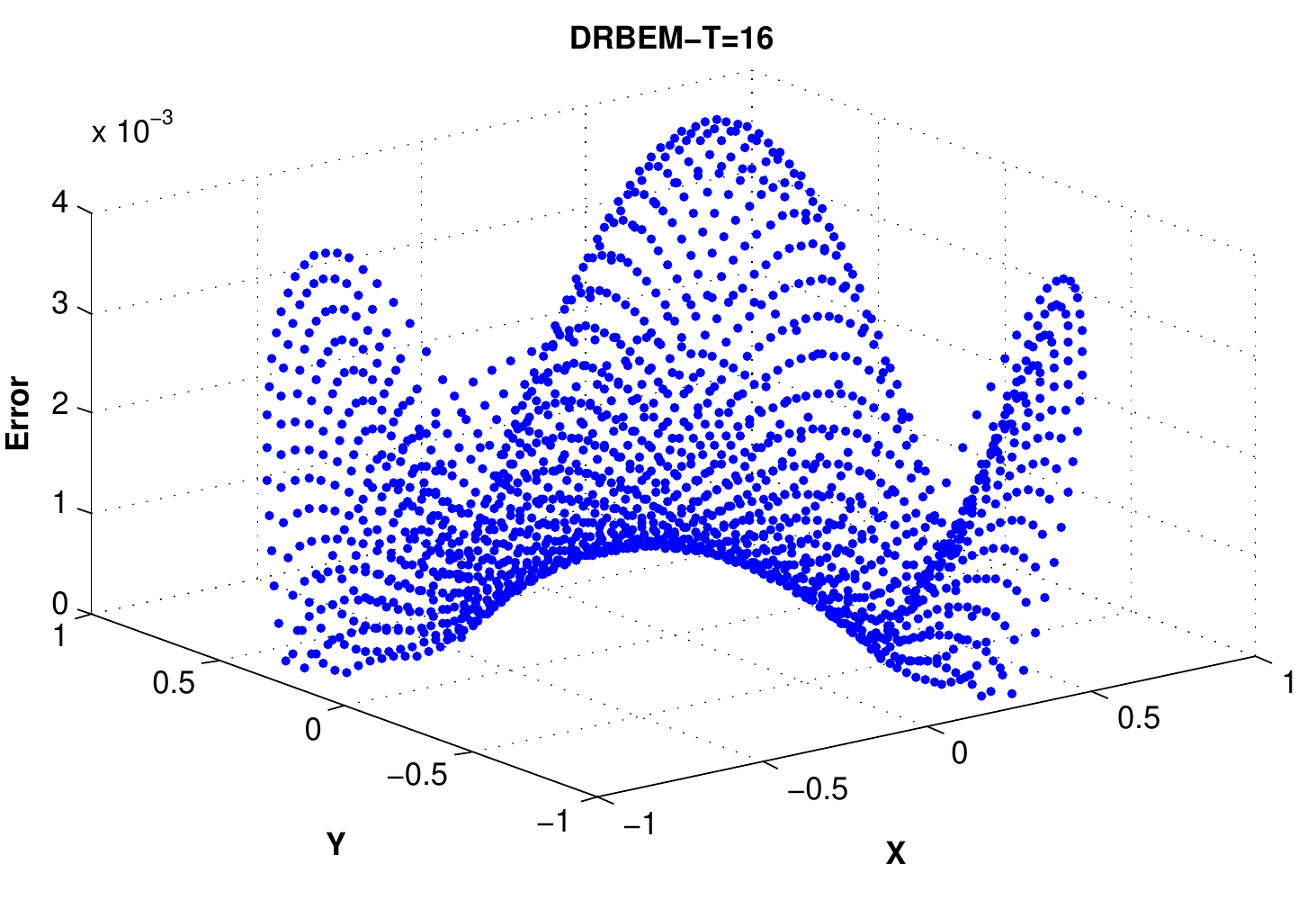}
\caption{Graphs of the average values of relative error with $\tau=1/100$, $N=750$ and $\alpha=1.95$ for Test problem 2.}
\end{figure}

\begin{figure}[t]
\includegraphics[width=8cm, height=7.cm]{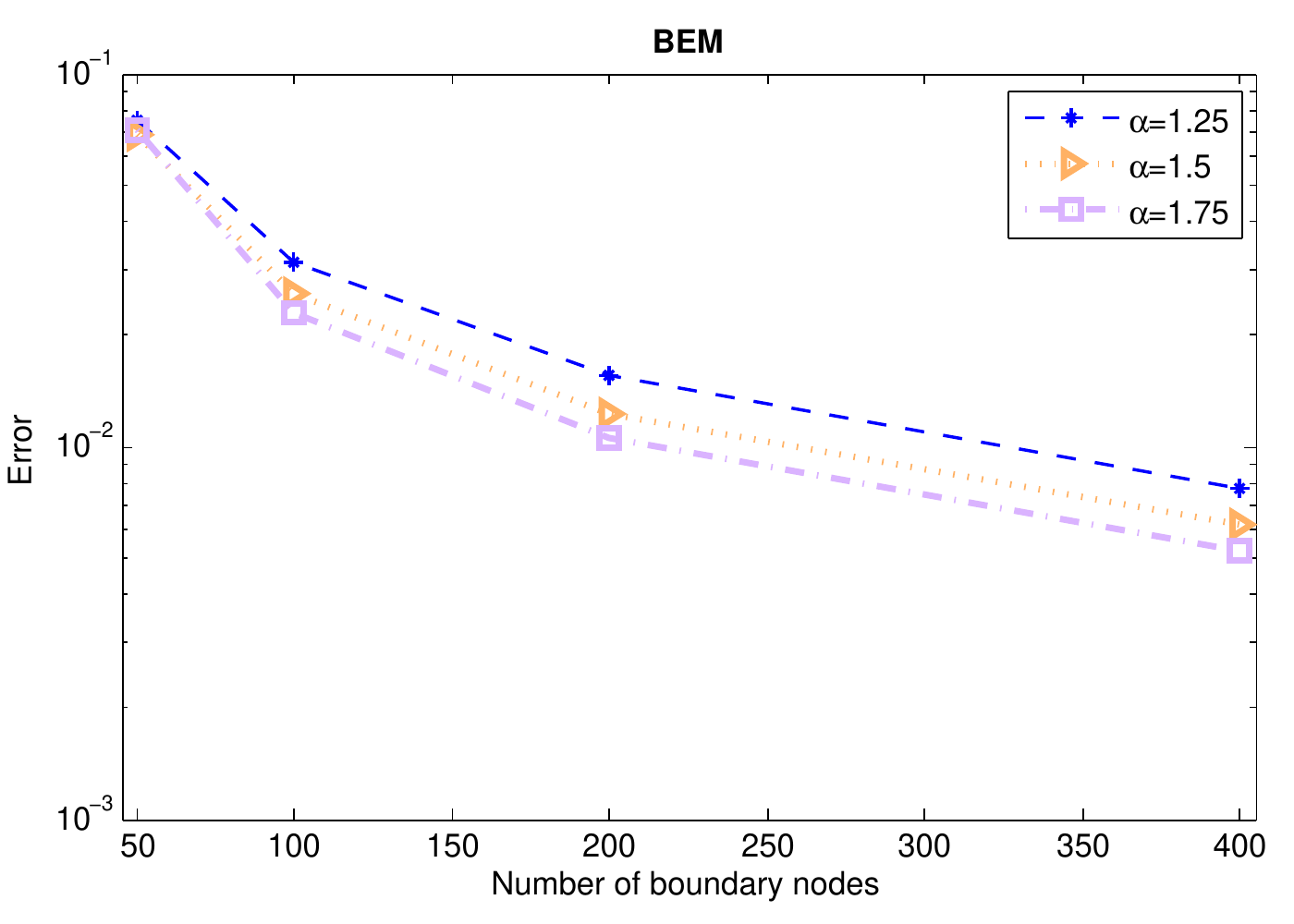}
\includegraphics[width=8cm, height=7.cm]{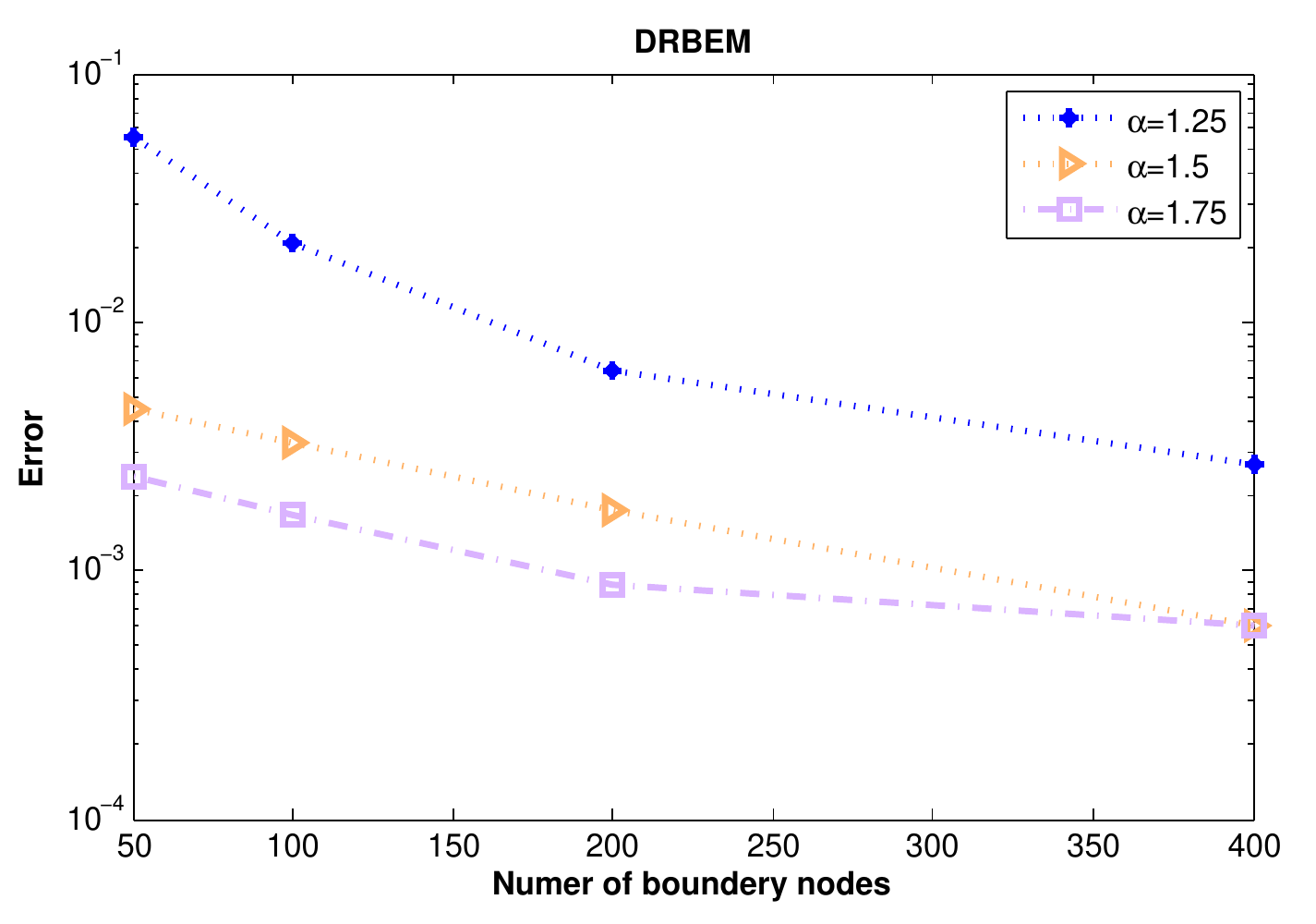}
\caption{Graphs the average values of relative errors with $\tau=1/10$ for Test problem 3.}
\end{figure}

\begin{figure}[t]
\includegraphics[width=8cm, height=7.cm]{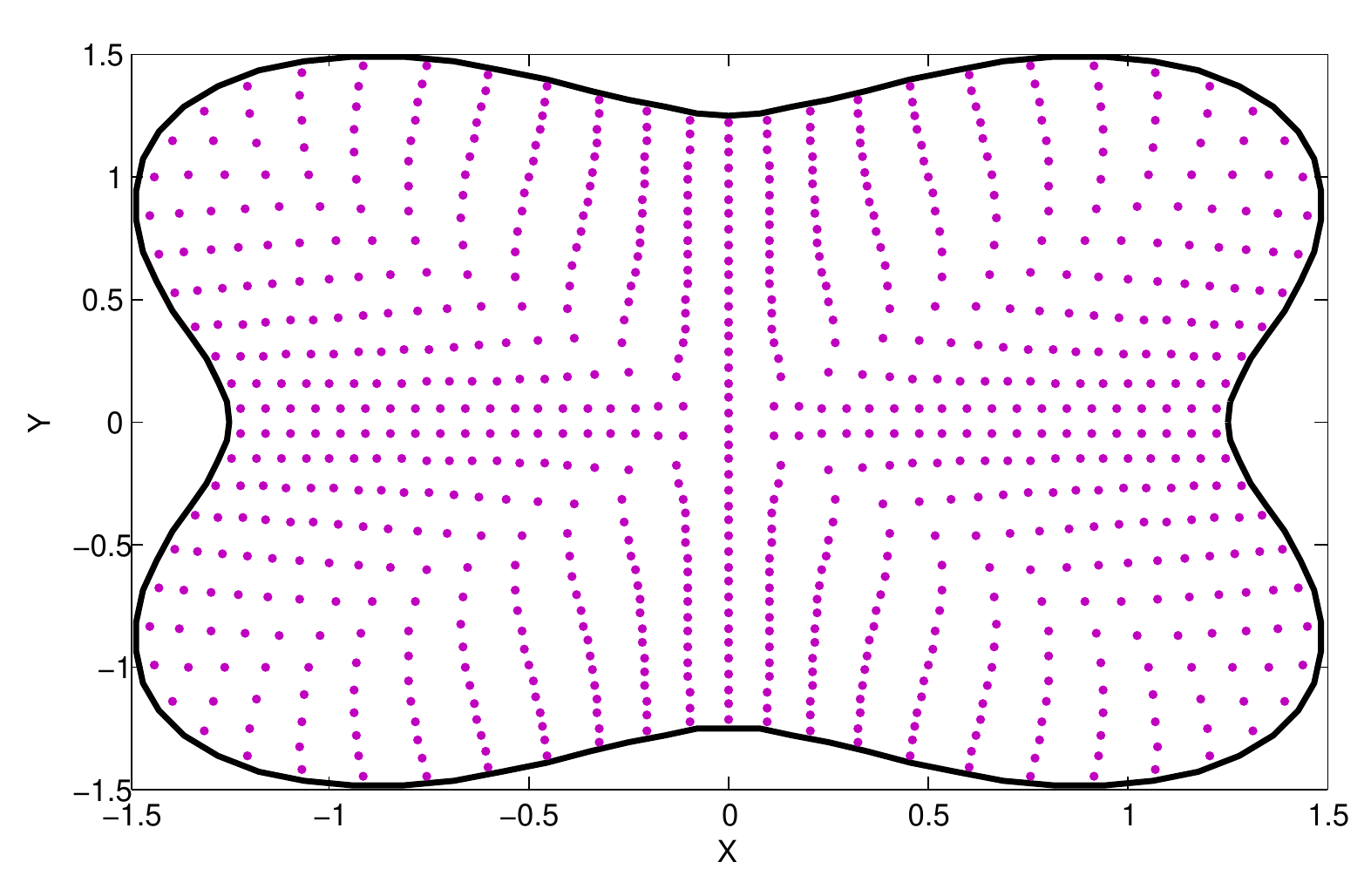}
\includegraphics[width=8cm, height=7.cm]{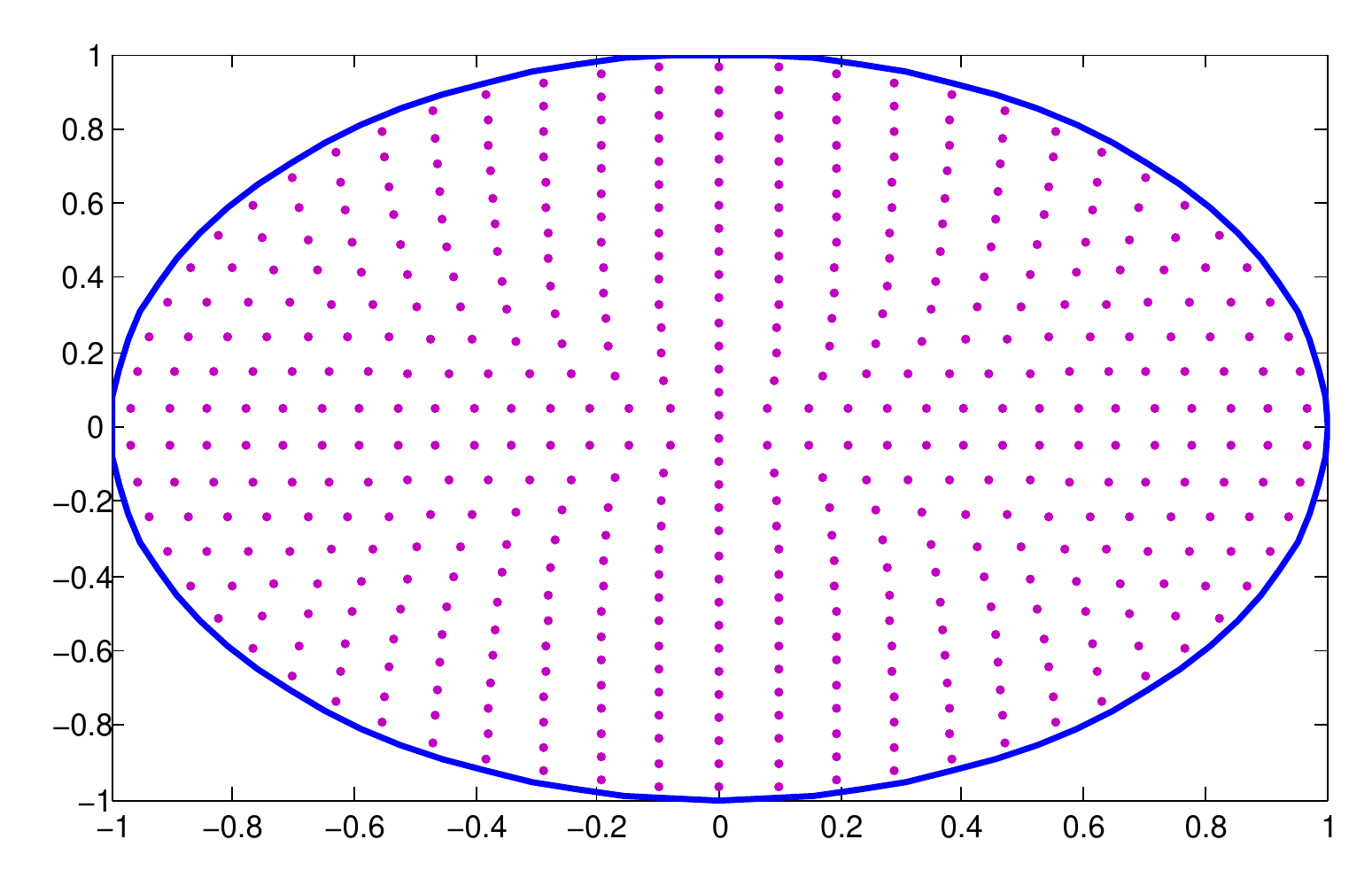}
\caption{Graphs the uniform mesh obtained for Test problem 3 (left plane) and Test problem 2 (right plane).}
\end{figure}


\begin{thebibliography}{99}
\bibitem{A. Mohebbi1} M. Abbaszadeh, A. Mohebbi, A fourth-order compact solution of the two-dimensional modified anomalous fractional sub-diffusion equation with a nonlinear source term, Comput. Math. Appl. 66 (2013) 1345-1359.

\bibitem{N. Alsoy-Akun}N. Alsoy-Akun, M. Tezer-Sezgin, DRBEM solution of the theorem-solutal buoyancy induced mixed convection flow problems, Eng. Anal. Bound. Elem. 37 (2013) 513-526.


\bibitem{Ang}W. T. Ang, K. C. Ang, A dual reciprocity boundary element solution of a generalized nonlinear Schr\"{o}dinger equation, Numer. Meth. Part. D. E. 20 (2004) 843-854.
\bibitem{Ang1}W. T. Ang, A Beginnes Course in Boundary Element Methods, Universal Publishers, Boca Raton, USA, 2007.


\bibitem{C. A. Brebbia}C. A. Brebbia, S. Walker, The Boundeary Element Techniques in Engineering, Newnes-Butterworths, London, 1980.

\bibitem{H. Brezis}H. Brezis, Functional Analysis, Sobolove Space and Partial Differential Equations, Springer New York, Dordrecht, Heidelberg, London, 2011.

\bibitem{C. Bozkaya}C. Bozkaya, Boundary element method solution of initial and boundary value problems in fluid dynamics, Ph.D Thesis, Technical University of Midell East, 2008.

\bibitem{C. Bozkaya1}C. Bozkaya, M. Tezer--Sezgin, A direct BEM solution to MHD flow in electrodynamically coupled rectangular channels, Comput. Fluids, 66 (2012) 177--182.

\bibitem{A. V}A. V. Chechkin, R. Gorenflo, I. M. Sokolov, Fractional diffusion in inhomogeneous media, J. Phys. A. 38 (2005) 679-684.



\bibitem{Cui1}M. Cui, Compact finite difference method for the fractional diffusion equation, J. Comput. Phys. 228 (2009) 7792-7804.

\bibitem{Cui2}M. Cui, A high-order compact exponential scheme for the fractional convection-diffusion equation, J. Comput. Appl. Math. 225 (2014) 404-416.

\bibitem{Dehghan_1} M. Dehghan, D. Mirzaei, The dual reciprocity boundary element method (DRBEM) for two-dimensional sine-Gordon equation, Comput. Methods Appl. Mech. Engrg. 197 (2008) 476-486.



\bibitem{K. Diethelm}K. Diethelm, N. J. Ford, Analysis of fractional differential equations, J. Math. Anal. Appl. 265 (2002) 229--248.


\bibitem{R. Du}R. Du, W. R. Cao, Z. Z. Sun, A compact difference scheme for the fractional diffusion-wave equation, Appl. Math. Model. 34 (2010) 2998-3007.

\bibitem{V. V. Gafiychuk}V. V. Gafiychuk, B. Yo. Datsko, Pattern formation in a fractional reaction diffusion systems, Physica A: Statistical Mechanics and its Applications, 365 (2006) 300-306.

\bibitem{R. Gorenflo}R. Gorenflo, F. Mainardi, D. Moretti, P. Paradis, Time fractional diffusion: A discrete random walk approach, Nonlinear Dyn. 29 (2002) 129-143.

\bibitem{Gumgum}S. G$\ddot{u}$mg$\ddot{u}$m, M. Tezer--Sezgin, DRBEM solution of mixed convection flow of nanofluids in enclosures with moving walls, J. Comput. Appl. Math. 259 (2014) 730--740.

\bibitem{H. Hosseinzadeh} H. Hosseinzadeh, M. Dehghan, D. Mirzaei, The boundary elements method for magnetohydrodynamic (MHD) channel flows at high Hartmann numbers, Appl. Math. Model. 37 (2013) 2337-2351.

\bibitem{H. Hosseinzadeh1} H. Hosseinzadeh, M. Dehghan, A new scheme based on boundary elements method to solve linear Helmholtz and semi--linear Poisson's equations, Eng. Anal. Bound. Elem. 43 (2014) 124--135.

\bibitem{H.Jiang}H. Jiang, F. Liu, I. Turner, K. Burrage, Analytical solutions for the multi--term time--space Caputo--Riesz fractional advection--diffusion equations on a finite domain, J. Math. Anal. Appl. 389 (2012) 1117-1127.

\bibitem{Y. Jiang}Y. Jiang, J. Ma, High-order finite element methods for time-fractional partial differential equations, J. Comput. Appl. Math. 235 (2011) 3285-3290.

\bibitem{J. Katsikadelis}J. T. Katsikadelis, Boundary Elements Methods, Theory and Application, Elsevier, 2002.

\bibitem{J. T. Katsikadelis}J. T. Katsikadelis, The analog equation method, A powerful BEM--based technique for solving linear and nonlinear engineering problems, boundary elements method, Southampton, Computational Mechanics, Publications, 1994.

\bibitem{Katsikadelis1}J. T. Katsikadelis, The BEM for numerical solution of partial fractional differential equations, Comput. Math. Appl. 62 (2011) 891--901.

\bibitem{Katsikadelis2}J. T. Katsikadelis, The fractional diffusion--wave equation in bounded inhomogeneous anisotropic media. An AEM Solution, Recent Advances in Boundary Element Methods 2009, pp 255--276.

\bibitem{P. K. Kythe}P. K. Kythe, Fundamental Solution for Differential Operators and Application, Birkhauser, Boston, 1996.

\bibitem{Liu}F. Liu, P. Zhuang, V. Anh, I. Turner, K. Burrage, Stability and convergence of the difference methods for the space--time fractional advection--diffusion equation, J. Appl. Math. Comput. 191 (2007) 12-20.

\bibitem{F. Liu}F. Liu, P. Zhuang, K. Burrage, Numerical methods and analysis for a class of fractional advection--dispersion models, Comput. Math. Appl. 64 (2012) 2990--3007.

\bibitem{F. Liu1}F. Liu, S. Shen, V. Anh, I. Turner, Analysis of a discrete non--Markovian random walk approximation for the time fractional diffusion equation, ANZIAM. J. 46 (2005) 488--504.

\bibitem{F. Mainardi3}F. Mainardi, The fundamental solutions for the fractional diffusion-wave equation, Appl. Math. Lett. 9 (1996) 23-28.

\bibitem{F. Mainardi4}F. Mainardi, The time fractional diffusion--wave equation, Radiophysics and Quantum Electronics, 38 (1995) 13--24.

\bibitem{Mark M. Meerschaert} M. M. Meerschaert, C. Tadjeran, Finite difference approximation for two--sided space--fractional partial differential equations, Appl. Numer. Math. 56 (2006) 80-90.

\bibitem{R. Metzler}R. Metzler, J. Klafter, The random walk’s guide to anomalous diffusion: a fractional dynamics approach, Phys. Rep. 339 (2000) 1-77.

\bibitem{D. Mirzaei} D. Mirzaei, M. Dehghan, Boundary element solution of the two--dimensional sine--Gordon equation using continuous linear elements, Eng. Anal. Bound. Elem. 34 (2010) 51--59.


\bibitem{Miller}K. S. Miller, B. Ross, An Introduction the Fractional Calculus and Fractional Differential Equations, New York and London, Academic Press, 1974.

\bibitem{A. Mohebbi}A. Mohebbi, M. Abbaszadeh, M. Dehghan, The use of a meshless techniqu based on collocation and radial basis functions for solving the fractional nonlinear Schrodinger equation arising in quantum mechanics, Eng. Anal. Bound. Elem. 37 (2013) 475-485.

\bibitem{A. Mohebbi2} A. Mohebbi, M. Abbaszadeh, M. Dehghan, Compact finite difference schemes and RBF meshless approach for solving 2D Rayleigh--Stokes problem for a heat generalized second grade fluid with fractional derivatives, Comput. Methods Appl. Mech. Engrg. 264 (2013) 163-177.

\bibitem{Oldham}K. B. Oldham, J. Spanier, The Fractional Calculus, New York and London, Academic Press, 1974.

\bibitem{I. Padlubny}I. Podlubny, Fractional Differential Equations, Academic Press, San Diego, 1999.


\bibitem{W. R. Schneider}W. R. Schneider, W. Wyss, Fractional diffusion and wave equations, J. Math. Phys. 30 (1989) 134-144.

\bibitem{Z. Sedaghatjoo}Z. Sedaghatjoo, H. Adibi, Calculation of domain integrals of two dimensional boundary element method, Eng. Anal. Bound. Elem. 36 (2012) 1917-1922.

\bibitem{Z. Sedaghatjoo1}Z. Sedaghatjoo, M. Dehghan, H. Hosseinzadeh, The use of continuous boundary elements in the boundary elements method for domain with non--smooth boundaries via finite difference approach, Comput. Math. Appl. 65 (2013) 983--995.

\bibitem{Tadjeran}C. Tadjeran, M. M. Meerschaert, H. P. Scheffler, A second-order accurate numerical approximation for the fractional diffusion equation, J. Comput. Phys. 213 (2006) 205-213.

\bibitem{M. Tezer--Sezgin}M. Tezer--Sezgin, S. H. Aydin, Solution of magnetohydrodynamics flow problems using the boundary element method, Eng. Anal. Bound. Elem. 30 (2006) 411-418.

\bibitem{M. Tezer--Sezgin1}M. Tezer--Sezgin, S. Dest, Boundary element method for MHD channel flow with orbitrary wall conductivity, Appl. Math. Model. 18 (1994) 429--436.


\bibitem{W. Wyss}W. Wyss, Fractional diffusion equation, J. Math. Phys. 27 (1986) 2782-2785.

\bibitem{Limei Li}Limei Li, Da Xu, Man Luo, Alternating direction implicit Galerkin finite element method for the two-dimensional fractional diffusion--wave equation, J. Comput. Phys. 15 (2013) 471-485.

\bibitem{S.B. Yuste}S. B. Yuste, L. Acedo, An explicit finite difference method and a new von Neumann-type stability analysis for fractional diffusion equations, SIAM. J. Numer. Anal. 42 (5) (2005) 1862-1874.

\bibitem{Ya-Nan Zhang}Ya-Nan Zhang, Z. Z. Sun, Xuan Zhao, Compact alternating direction implicit for the two--dimensional diffusion--wave equation, SIAM, J. Numer. Anal. 50 (2012) 1535-1555.

\bibitem{Yang Zhang}Yang Zhang, A finite difference for fractional partial differential equation, Appl. Math. Comput. 215 (2009) 524-529.

\bibitem{Zhuang}P. Zhuang, F. Liu, V. Anh, I. Turner, New solution and analytical techniques of the implicit numerical methods for the anomalous sub--diffusion equation, SIAM J. Numer. Anal. 46  (2008) 1079-1095.

\bibitem{Zhi-zhong Sun}Z. Z. Sun, Xiaonan Wu, A fully discrete difference for a diffusion-wave system, Appl. Numer. Math. 56 (2006) 193-209.


\end{thebibliography}
\end{document}